\documentclass{lms}

\usepackage{amssymb,amsmath} 
\usepackage[mathscr]{euscript} 
\usepackage[all]{xy}

\newcommand{\C}{\ensuremath{\mathbb{C}}}
\newcommand{\N}{\ensuremath{\mathbb{N}}}

\renewcommand{\O}{\mathcal{O}}
 
\renewcommand{\Im}{\mathop{\rm Im}}
 \newcommand{\Ker}{\mathop{\rm Ker}}
\newcommand{\iN}{\hbox{ {\leaders\hrule\hskip.2cm}{\vrule height .22cm} }}
\newcommand{\pa}[3][]{ \frac{ \partial^{#1} {#2} }{ \partial {#3}^{#1} } }
 
\newcommand{\at}[1]{\big|_{#1}}

\newcommand{\gC}{\ensuremath{\mathfrak{g}^{\C}}}
\renewcommand{\u}{\ensuremath{\mathfrak{u}}}
\newcommand{\su}{\ensuremath{\mathfrak{su}}}
\renewcommand{\sl}{\ensuremath{\mathfrak{sl}}} \renewcommand{\phi}{\varphi}
 \newcommand{\ad}{\mathop{\rm ad}}

\newtheorem{theorem}{Theorem}[section]
\newtheorem{proposition}[theorem]{Proposition}
\newtheorem{lemma}[theorem]{Lemma} 
\newtheorem{corollary}[theorem]{Corollary}

\newnumbered{definition}[theorem]{Definition}
\newnumbered{example}[theorem]{Example}
\newunnumbered{remark}[theorem]{Remark}

\title{Hamiltonian stationary tori in the complex projective plane}
\author{Fr\'ed\'eric H\'elein \and Pascal Romon}

\classno{53C55 (primary), 53C42, 53C25, 58E12 (secondary)}

\begin{document}
\maketitle

\section{Introduction}
Hamiltonian stationary Lagrangian surfaces are Lagrangian surfaces of a given
four-dimensional manifold endowed with a symplectic and a Riemannian structure,
which are critical points of the area functional with respect to a particular
class of infinitesimal variations preserving the Lagrangian constraint: 
the compactly supported Hamiltonian vector fields.  The
Euler--Lagrange equations of this variational problem are highly simplified when
we assume that the ambient manifold ${\cal N}$ is K\"ahler.  In that case we can
make sense of a {\em Lagrangian angle function} $\beta$ along any
simply-connected Lagrangian submanifold $\Sigma \subset {\cal N}$ (uniquely
defined up to the addition of a constant).  And as shown in \cite{ScWo} the mean
curvature vector of the submanifold is then $\vec{H}= J\,\nabla \beta$, where
$J$ is the complex structure on ${\cal N}$ and $\nabla \beta$ is the gradient of
$\beta$ along $\Sigma$. It turns out that $\Sigma$ is Hamiltonian stationary if and
only if $\beta$ is a harmonic function on $\Sigma$.\\

\noindent A particular subclass of solutions occurs when $\beta$ is constant:
the Lagrangian submanifold is then simply a minimal one.  In the case where
${\cal N}$ is a Calabi--Aubin--Yau manifold, such submanifolds admit an
alternative characterization as {\em special Lagrangian}, a notion which has
been extensively studied recently because of its connection with string theories
and the mirror conjecture, see \cite{SYZ}.\\

\noindent An analytical theory of two-dimensional Hamiltonian stationary
Lagrangian submanifolds was constructed by R. Schoen and J. Wolfson \cite{ScWo},
proving the existence and the partial regularity of minimizers.  In contrast our
results in the present paper rest on the fact that, for particular ambient
manifolds ${\cal N}$, Hamiltonian stationary Lagrangian surfaces are solutions
of an integrable system.  This was discovered first in the case when ${\cal
N}=\C^2$ in \cite{HR1} and \cite{HR2}.  In a subsequent paper \cite{HR3} we
proved that the same problem is also completely integrable if we replace $\C^2$
by any two-dimensional Hermitian symmetric space.  Among these symmetric spaces
one very interesting example is $\C P^2$, because any simply-connected
Lagrangian surface in $\C P^2$ can be lifted into a Legendrian surface in $S^5$.
Furthermore the cone in $\C^3$ over this Legendrian surface is actually a
singular Lagrangian three-dimensional submanifold in $\C^3$; and the cone in
$\C^3$ is Hamiltonian stationary if and only if the surface in $\C P^2$ is so.\\

\noindent A similar correspondence has been remarked and used in \cite{J},
\cite{McI3} and \cite{Has} in the case of minimal Lagrangian surfaces in $\C
P^2$ and allows these Authors to connect results on minimal Lagrangian surfaces
in $\C P^2$ \cite{Sh} to minimal Legendrian surfaces in $S^5$ \cite{MM} and
special Lagrangian cones in $\C^3$.\\

\noindent Our aim in this paper is the following:
\begin{itemize}
\item to expound in details the correspondence between Hamiltonian stationary
Lagrangian surfaces in $\C P^2$ and Hamiltonian stationary Legendrian surfaces
in $S^5$ and a formulation using a family of curvature free connections of this
integrable system (theorem \ref{1.6.theobis}).  
We revisit here the formulation given in \cite{HR3}, using {\em twisted loop groups}.  
Roughly speaking it rests on the identifications 
$\C P^2\simeq SU(3)/S(U(2)\times U(1))$ and $(S^5,\hbox{ contact structure}) \simeq
\left( U(3)/U(2)\times U(1), A^3_3=0\right)$, where $A^3_3$ is a component of
the Maurer--Cartan form.  We also show that this problem has an alternative
formulation, analogous to the theory of K. Uhlenbeck \cite{U} for harmonic maps
into $U(n)$, using {\em based loop groups}.

\item to define the notion of {\em finite type} Hamiltonian stationary
Legendrian surfaces in $S^5$: we give here again two definitions, in terms of
twisted loop groups (which is an analogue to the description of finite type
harmonic maps into homogeneous manifolds according to \cite{BP}) and in terms of
based loop groups (an analogue to the description of finite type harmonic maps
into Lie groups according to \cite{BFPP}).  We prove the equivalence between the
two definitions because we actually need this result for the following.  We
believe that this fact should be well known to some specialists in the harmonic
maps theory, but we did not find it in the literature.

\item we prove in theorem \ref{4.MainTheo} that all Hamiltonian stationary 
Lagrangian tori in $\C P^2$ (and
hence Hamiltonian stationary Legendrian tori in $S^5$) are of finite type.  This
is the main result of this paper.  Our proof focuses on the case of Hamiltonian
stationary tori which are not minimal, since the minimal case has been studied
by many authors (\cite{BFPP}, \cite{Sh}, \cite{MM}, \cite{Has},
\cite{McI1},\cite{McI2}, \cite{McI3}, \cite{J}).  The method here is adapted
from the similar result for harmonic maps into Lie groups in \cite{BFPP}.
However the strategy differs slightly: we use actually the two existing
formulations of finite type solutions, using twisted or based loop groups.  One
crucial step indeed is the construction of a {\em formal Killing field},
starting from a given torus solution.  This step can be slightly simplified here
in the twisted loop groups formulation, because the semi-simple element we start
with is then just constant.  However proving that the formal Killing field is
{\em adapted} requires more work in the twisted loop groups formulation
(actually we were not able to do it directly) than in the based loop groups
formulation; here we take advantage from the two formulations to avoid the
difficulties and to conclude.

\item lastly we give some examples of Hamiltonian stationary Legendrian tori in
$S^5$: we construct in theorem \ref{homogeneous} a family of solutions which are equivariant in some sense
under the action of the torus, that we call {\em homogeneous} Hamiltonian
stationary tori.  These are the simplest examples that one can build.
\end{itemize}

\noindent Let us add that the structure of the integrable system studied here
fits in a classification of elliptic integrable systems proposed by C.L. Terng
\cite{T}, as a 2nd $(U(3),\sigma,\tau)$-system\footnote{We have here exchanged
the notations $\sigma$ and $\tau$ with respect to \cite{T} in order to be
consistent with our notations in \cite{HR3}.}, where $\sigma$ is an involution
of $U(3)$ such that its fixed set is $U(3)^\sigma\simeq U(2)\times U(1)$ and
$U(3)/U(2)^\sigma\simeq \C P^2$ and $\tau$ is a 4th order automorphism (actually
$\tau^2=\sigma$) which encodes the symplectic structure on $\C P^2$ or the
Legendrian structure on $S^5$.\\

\noindent {\em Notations} --- For any matrix $M\in GL(n,\C)$, we denote by
$M^\dagger :=\,^t\overline{M}$.

\section{Geometrical description of Hamiltonian stationary Lagrangian surfaces
in $\mathbb{C}P^n$}

\subsection{The Lagrangian angle}

The complex projective space $\mathbb{C}P^n$ can be identified with the quotient
manifold $S^{2n+1}/S^1$. It is a complex manifold with complex structure $J$. We
denote by $\pi:S^{2n+1}\longrightarrow \mathbb{C}P^n$
 the canonical projection a.k.a. Hopf fibration, and equip $\mathbb{C}P^n$ with the
Fubini-Study Hermitian metric, denoted by $\langle \cdot , \cdot \rangle _{\mathbb{C}P^n} =
\langle \cdot , \cdot \rangle  - i\omega(\cdot , \cdot )$, where
$\langle \cdot , \cdot \rangle$ is a Riemannian metric and $\omega$ is the
K\"ahler form\footnote{note that the sign
convention may vary in the literature, e.g.\,some Authors use $\langle \cdot ,
\cdot \rangle _{\mathbb{C}P^n} = \langle \cdot , \cdot \rangle +
i\omega(\cdot , \cdot )$.}. 
For each $z\in S^{2n+1}$ we let ${\cal H}_z$ be the complex $n$-subspace 
in $T_{z} S^{2 n+1}\subset \mathbb{C}^{n+1}$ which is Hermitian orthogonal to $z$ 
(and hence to the fiber of $d\pi_{z}$). By construction of the Fubini-Study metric,
$d\pi_z:{\cal H}_z\longrightarrow T_{\pi(z)}\mathbb{C}P^n$ is 
an isometry between complex Hermitian spaces. We call the
subbundle ${\cal H}:= \cup_{z\in S^{2n+1}}{\cal H}_z$ of $T S^{2n+1}$
the {\em horizontal distribution}. It defines in a natural way a connection
$\nabla^\mathit{Hopf}\simeq \nabla ^H$ on the Hopf bundle 
$\pi:S^{2n+1}\longrightarrow \mathbb{C}P^n$, whose curvature is $2i\omega$. 
As a consequence \cite{Re2}:
\begin{proposition}\label{2.1.prop}
Let $\Omega$ be a {\em simply connected} open subset of $\mathbb{R}^n$ and 
$u:\Omega \longrightarrow \mathbb{C}P^n$
be a smooth Lagrangian immersion, i.e.\,such that $u^*\omega = 0$. Then there exists a lift
\begin{center}
\hskip2cm
\xymatrix{
& S^{2n+1} \ar@{->}[d]^\pi \\
\Omega \ar@{->}[ur]^{\widehat{u}} \ar@{->}[r]^u & \mathbb{C}P^n
}
\end{center}
such that $\left(u^*\nabla^H\right) \widehat{u} = 0$ (where $u^*\nabla^H$ is 
the pull-back by $u$ of the connection $\nabla^H$). This lift is unique up to multiplication
by a unit complex number. Moreover the pull-back by $\widehat{u}$ of the symplectic
form  $\omega$ on $\mathbb{C}^{n+1}$ vanishes; we say that $\hat{u}$ is 
\emph{Legendrian}.
\end{proposition}

Taking $u,\hat{u}$ as above, we define, for any orthonormal framing
$(e_{1},\ldots,e_{n})$ of $T\Omega$, the Lagrangian angle $\beta$ by
\[
e^{i \beta} = dz^1 \wedge \ldots \wedge dz^{n+1} \big( \hat{u}, d\hat{u}(e_{1}),
\ldots, d\hat{u}(e_{n}) \big) .
\]
which makes sense because $( \hat{u}, d\hat{u}(e_{1}),\ldots, d\hat{u}(e_{n}))$ is
a Hermitian-orthonormal frame, for any $x\in\Omega$. Furthermore, the result is
independent from the choice of the framing, and depends on the choice of the lift
$\hat{u}$ only through multiplication by a unit complex constant. Hence $\beta$ is
defined up to an additive constant and $d\beta$ is always well-defined along any
Lagrangian immersion $u$. Another characteristic property of the Lagrangian angle
relates it to the mean curvature vector field $\vec{H}$ along $u$:
\begin{equation}\label{2.2.H=JDb}
\vec{H} = {1\over n}J\nabla \beta
\end{equation}
or equivalently $d\beta = - \vec{H} \iN \omega$ (see~\cite{Br,Da} for details). 

\subsection{Hamiltonian stationary Lagrangian submanifolds}
A \emph{Hamiltonian stationary} Lagrangian submanifold $\Sigma$ in
$\mathbb{C}P^n$ is a Lagrangian submanifold which is a critical point of the
$n$-volume functional ${\cal A}$ under first variations which are {\em
Hamiltonian vector fields} with compact support.  This means that for any smooth
function with compact support $h\in {\cal
C}^\infty_c(\mathbb{C}P^n,\mathbb{R})$, we have
\[
\delta {\cal A}_{\xi_h}(\Sigma) := \int_\Sigma \left\langle \vec{H},
\xi_h\right\rangle_E d\hbox{vol} = 0,
\]
where $\xi_h$ is the Hamiltonian vector field of $h$, i.e.\,satisfies
$\xi_h \iN \omega +dh = 0$ or $\xi_h = J\nabla h$.  We also remark that if
$f\in {\cal C}^\infty_c(\Sigma,\mathbb{R})$, then there exist smooth extensions
with compact support $h$ of $f$, i.e.\,functions $h\in {\cal
C}^\infty_c(\mathbb{C}P^n,\mathbb{R})$ such that $h_{|\Sigma} = f$, and moreover
the normal component of $\left(\xi_h\right)_{|\Sigma}$ does not depend on the
choice of the extension $h$ (it coincides actually with $J\nabla f$, where
$\nabla$ is here the gradient with respect to the induced metric on $\Sigma$).
So we deduce from above that $\delta {\cal
A}_{\xi_h}(\Sigma) = {1\over n}\int_\Sigma \left\langle \nabla \beta, \nabla
f\right\rangle_E d\hbox{vol}$.  This implies the following.
\begin{corollary}
Any Lagrangian submanifold $\Sigma$ in $\mathbb{C}P^n$ is Hamiltonian stationary
if and only if $\beta$ is a harmonic function on $\Sigma$, i.e.
\[
\Delta_\Sigma \beta = 0.
\]
\end{corollary}

\noindent This theory extends to non simply connected surfaces $\Sigma$
with the following restrictions.  Let $\gamma$ be a homotopically non trivial
loop.  The Legendrian lift of $\gamma$ needs not close, so that in general its
endpoints $p_1, p_2 \in S^{2n+1}$ are multiples of each other by a factor $e^{i
\theta}$.  The same holds for the Lagrangian angle: $\beta ( p_2 ) \equiv \beta
( p_1 ) + (n+1) \theta \bmod 2 \pi$ (since the tangent plane is also shifted by
the Decktransformation $z \longmapsto e^{i \theta} z$).  In particular $\beta$
is not always globally defined on surfaces in $\mathbb{C} P^n$ with non trivial
topology, unless the Legendrian lift is globally defined in $S^{2n+1} /
\mathbb{Z}_{n+1}$ (here $\mathbb{Z}_{n+1}$ stands for the $n+1$-st roots of
unity in $SU(n+1)$).

\subsection{Conformal Lagrangian immersions into $\mathbb{C}P^2$}
We now set $n= 2$.  We suppose that $\Omega$ is a simply connected open subset
of $\mathbb{R}^2\simeq \mathbb{C}$ and consider a conformal Lagrangian immersion
$u:\Omega\longrightarrow \mathbb{C}P^2$.  This implies that we can find a
function $\rho:\Omega\longrightarrow \mathbb{R}$ and two sections $E_1$ and
$E_2$ of $u^*T\mathbb{C}P^2$ such that $\forall (x,y)\in \Omega$,
$(E_1(x,y),E_2(x,y))$ is an Euclidean orthogonal basis over $\mathbb{R}$ of
$T_{u(x,y)}u(\Omega)$ and
\[
d u = e^\rho\left( E_1d x + E_2d y\right).
\]
We observe that, due to the fact that $u$ is Lagrangian, $(E_1,E_2)$ is a also a
{\em Hermitian} basis over $\mathbb{C}$ of $T_{u(x,y)}\mathbb{C}P^2$.\\

\noindent Let \xymatrix{ \Omega\ar@{->}[r]^{\widehat{u}}\ar@{->}[dr]^u & S^5
\ar@{->}[d]^\pi \\
 & \mathbb{C}P^2 } be a parallel lift of $u$ as in Proposition \ref{2.1.prop}
 and $(e_1,e_2)$ be the unique section of $\widehat{u}^*{\cal H}\times
 \widehat{u}^*{\cal H}$ which lifts\footnote{recall that the condition that
 $v\in \left(\widehat{u}^*{\cal H}\right)_{(x,y)}$ means that $v$ is in the
 horizontal subspace ${\cal H}_{\widehat{u}(x,y)}$} $(E_1,E_2)$.  Then we have
\begin{equation}\label{duchapeau}
d\widehat{u} = e^\rho\left( e_1d x + e_2dy\right).
\end{equation}
Note that $\forall (x,y)\in \Omega$, $(e_1(x,y),e_2(x,y))$ is a Hermitian basis
of ${\cal H}_{\widehat{u}(x,y)}$, which is Hermitian orthogonal to
$\widehat{u}(x,y)$.  Hence $\forall (x,y)\in \Omega$,
$(e_1(x,y),e_2(x,y),\widehat{u}(x,y))$ is a Hermitian basis of $\mathbb{C}^3$.
Thus this triplet can be identified with some $\widehat{F}(x,y)\in U(3)$.  We
hence get the diagram \xymatrix{ & U(3)\ar@{->}[d]^{(\cdot \cdot *)}\\
\Omega \ar@{->}[ur]^{\widehat{F}} \ar@{->}[r]^{\widehat{u}} \ar@{->}[dr]^u &
S^5\ar@{->}[d]^{\pi}\\
& \mathbb{C}P^2 }, where $(\cdot \cdot *)$ is the mapping
$(e_1,e_2,e_3)\longmapsto e_3$.\\

\noindent We define the Maurer--Cartan form $\widehat{A}$ to be the 1-form on
$\Omega$ with coefficients in $\u(3)$ such that $d\widehat{F} = \widehat{F}\cdot
\widehat{A}$.  Then we remark that the horizontality assumption $\langle
d\widehat{u},\widehat{u}\rangle _{\mathbb{C}^3} = 0$ exactly means that
\begin{equation}\label{A33=0}
\widehat{A}^3_3 = 0.
\end{equation}
Moreover the Lagrangian angle function $\beta_{\widehat{u}}$ along $\widehat{u}$
can be computed by
\[
e^{i\beta_{\widehat{u}}} = dz^1 \wedge dz^2 \wedge dz^3 (e_1,e_2,\widehat{u}) 
= \det \widehat{F}.
\]
As in \cite{HR1} we consider a larger class of framings of $u$ as follows.
\begin{definition}
A \emph{Legendrian framing of $u$ along $\widehat{u}$} is a map
$F:\Omega\longrightarrow U(3)$ such that
\begin{itemize}
\item $(\cdot \cdot *)\circ F = \widehat{u}$ \item $\det F =
e^{i\beta_{\widehat{u}}}$.
\end{itemize}
\end{definition}
It is easily seen that the first condition is equivalent to the fact that there
exists a smooth map $G:\Omega \longrightarrow U(3)$ (a gauge transformation) of
the type
\[
G(x,y) = \left(\begin{array}{cc}g(x,y) & 0 \\ 0 & 1\end{array}\right),\quad
\hbox{where }g:\Omega\longrightarrow U(2)
\]
such that \[F(x,y) = \widehat{F}(x,y)\cdot G^{-1}(x,y).
\]
And then the second one is equivalent to say that $g$ takes values in $SU(2)$.

\subsection{A splitting of the Maurer--Cartan form of a Legendrian framing}
Using (\ref{duchapeau}) and (\ref{A33=0}) one obtains the following
decomposition of $\widehat{A}$:
\[
\widehat{A} = \widehat{A}_{\u(1)} + \widehat{A}_{\su(2)} +
\widehat{A}_{\mathbb{C}^2},
\]
with the notations
\[
\widehat{A}_{\u(1)} = \left(\begin{array}{cc}\widehat{\alpha}_{\u(1)} & 0\\0 &
0\end{array}\right),\quad \widehat{A}_{\su(2)} =
\left(\begin{array}{cc}\widehat{\alpha}_{\su(2)} & 0\\0 & 0\end{array}\right),
\]
\[
\hbox{and}\quad \widehat{A}_{\mathbb{C}^2} = e^\rho\left(\begin{array}{cc} 0 &
\epsilon dz + \overline{\epsilon}d\bar{z}\\
-\,^t\!\left(\epsilon dz + \overline{\epsilon}d\bar{z}\right) &
0\end{array}\right),
\]
where $\widehat{\alpha}_{\u(1)}$ is a 1-form on $\Omega$ with coefficients in
$\u(1)\simeq \mathbb{R}\left(\begin{array}{cc}i&0\\0&i\end{array}\right)$,
$\widehat{\alpha}_{\su(2)}$ a 1-form on $\Omega$ with coefficients in $\su(2)$,
$\epsilon:={1\over 2}\left(\begin{array}{c}1\\-i\end{array}\right)$ and
$\overline{\epsilon}:={1\over 2}\left(\begin{array}{c}1\\i\end{array}\right)$
(so that $\epsilon dz +
\overline{\epsilon}d\bar{z}=\left(\begin{array}{c}dx\\dy\end{array}\right)$).
Note that $\det F = e^{i\beta_{\widehat{u}}}$ implies $\widehat{\alpha}_{\u(1)}
= {d\beta_{\widehat{u}}\over
2}\left(\begin{array}{cc}i&0\\0&i\end{array}\right)$.\\
\noindent Now we let $F:\Omega\longrightarrow U(3)$ be a Legendrian framing and
$A:= F^{-1}\cdot dF$.  The relation $F = \widehat{F}\cdot G^{-1}$ implies that
$A = G\cdot \widehat{A}\cdot G^{-1} - dG\cdot G^{-1}$.  Hence
\[
A = A_{\u(1)} + A_{\su(2)} + A_{\mathbb{C}^2},
\]
where, using the fact that $\u(1)$ commutes with $\su(2)$,
\[
A_{\u(1)} = \left(\begin{array}{cc}\alpha_{\u(1)} & 0\\0 & 0\end{array}\right) =
\left(\begin{array}{cc}\widehat{\alpha}_{\u(1)} & 0\\0 & 0\end{array}\right) =
\left(\begin{array}{ccc}i{d\beta_{\widehat{u}}\over
2}&0&0\\0&i{d\beta_{\widehat{u}}\over 2}&0\\
0&0&0\end{array}\right),
\]
\[
A_{\su(2)} = \left(\begin{array}{cc}\alpha_{\su(2)} & 0\\0 & 0\end{array}\right)
= \left(\begin{array}{cc}g\cdot \widehat{\alpha}_{\su(2)}\cdot g^{-1} -dg\cdot
g^{-1} & 0\\0 & 0\end{array}\right)
\]
and
\[
A_{\mathbb{C}^2} = e^\rho\left(\begin{array}{cc} 0 & g\cdot\left( \epsilon dz +
\overline{\epsilon}d\bar{z}\right)\\
-\left(g\cdot \left(\epsilon dz +
\overline{\epsilon}d\bar{z}\right)\right)^\dagger & 0\end{array}\right).
\]
We can further split the last term $A_{\mathbb{C}^2}$ along $dz$ and $d\bar{z}$
as $A_{\mathbb{C}^2} = A_{\mathbb{C}^2}' + A_{\mathbb{C}^2}''$ where
\[
A_{\mathbb{C}^2}':= e^\rho\left(\begin{array}{cc} 0 & g\cdot \epsilon \\
-\left(g\cdot \overline{\epsilon}\right)^\dagger & 0\end{array}\right) dz
\quad\hbox{and}\quad A_{\mathbb{C}^2}'':= e^\rho\left(\begin{array}{cc} 0 &
g\cdot \overline{\epsilon} \\
-\left(g\cdot \epsilon\right)^\dagger & 0\end{array}\right) d\bar{z}.
\]

\subsection{Interpretation in terms of an automorphism}
As expounded in \cite{HR1} and \cite{HR3} the key point in order to exploit the
structure of an integrable system is to observe that the splitting $A =
A_{\u(1)} + A_{\su(2)} + A_{\mathbb{C}^2}' + A_{\mathbb{C}^2}''$ corresponds to
a decomposition along the eigenspaces of the following automorphism in
$\u(3)^{\mathbb{C}}$, the complexification\footnote{we can define
$\u(3)^{\mathbb{C}}$ as the set $M(3,\mathbb{C})$ with its standard complex
structure and with the conjugation mapping $c:M\longmapsto -M^\dagger$; clearly
$c$ is a Lie algebra automorphism, an involution and the set of fixed points of
$c$ is $\u(3)$.\\
Similarly the complexification $U(3)^{\mathbb{C}}$ is the set $GL(3,\mathbb{C})$
with its standard complex structure and the conjugation map $C:G\longmapsto
\left(G^\dagger\right)^{-1}$.} of $\u(3)$.  We let $J:=
\left(\begin{array}{cc}0&-1\\1&0\end{array}\right)$ and
\[
\begin{array}{cccl}
\tau: & \u(3)^{\mathbb{C}} & \longrightarrow & \u(3)^{\mathbb{C}}\\
 & M & \longmapsto & - \left( \begin{array}{cc}-J & 0\\0 & 1\end{array}\right)
 \cdot \,^t\!M\cdot \left( \begin{array}{cc}J & 0\\0 & 1\end{array}\right).
\end{array}
\]
It is then straightforward that $\tau$ is a Lie algebra automorphism, that
$\u(3)$ is stable by $\tau$ and that $\tau^4 = Id$.  Hence we can diagonalize
the action of $\tau$ over $\u(3)^{\mathbb{C}}$ and in the following we denote by
$\u(3)^{\mathbb{C}}_a$ the eigenspace of $\tau$ for the eigenvalue $i^a$, for $a
= -1,0,1,2$.  We first point out that the eigenspaces $\u(3)^{\mathbb{C}}_0$ and
$\u(3)^{\mathbb{C}}_2$, with eigenvalues 1 and $-1$ respectively, are the
complexifications of $\u(3)_0$ and $\u(3)_2$ respectively, where
\[
\u(3)_0:= \left\{ \left(\begin{array}{cc}g & 0 \\0 & 0\end{array}\right) /g\in
\su(2)\right\} \hbox{ and } \u(3)_2:= \left\{ \left(\begin{array}{ccc}\lambda
i&0&0\\0&\lambda i&0\\0&0&\mu i\end{array}\right) / \lambda, \mu\in
\mathbb{R}\right\}.
\]
This can be obtained by first computing that
\[
\tau \left(\begin{array}{cc}A & X \\-\,^t\!Y & d\end{array}\right) =
\left(\begin{array}{cc}J\,^t\!AJ & -JY \\-\,^t\!XJ & -d\end{array}\right), \quad
\forall A\in M(2,\mathbb{C}),\forall X,Y\in \mathbb{C}^2,\forall d\in
\mathbb{C},
\]
and by using the fact that $\forall A\in \sl(2,\mathbb{C})$, $J\,^t\!AJ = A$.
Similarly the eigenspaces $\u(3)^{\mathbb{C}}_{1}$ and
$\u(3)^{\mathbb{C}}_{-1}$, with eigenvalues $i$ and $-i$ respectively, are found
to be
\[
\u(3)^{\mathbb{C}}_{1} = \left\{ \left(\begin{array}{cc}0 & X \\-\,^t\!Y &
0\end{array}\right) /X,Y\in \mathbb{C}^2, JY = -iX\right\}
\]
and
\[
\u(3)^{\mathbb{C}}_{-1} = \left\{ \left(\begin{array}{cc}0 & X \\-\,^t\!Y &
0\end{array}\right) /X,Y\in \mathbb{C}^2, JY = iX\right\}.
\]
Now we have the following
\begin{lemma}
The eigenspaces $\u(3)^{\mathbb{C}}_{1}$ and $\u(3)^{\mathbb{C}}_{-1}$ can be
characterized by
\[
\u(3)^{\mathbb{C}}_{1} = \left\{\lambda \left(\begin{array}{cc}0 & h\cdot
\overline{\epsilon} \\
-\left(h\cdot \epsilon\right)^\dagger & 0\end{array}\right) / \lambda\in
[0,\infty), h\in SU(2)\right\}
\]
and
\[
\u(3)^{\mathbb{C}}_{-1} = \left\{ \lambda\left(\begin{array}{cc}0 & h\cdot
\epsilon \\
-\left(h\cdot \overline{\epsilon}\right)^\dagger & 0\end{array}\right) /
\lambda\in [0,\infty), h\in SU(2)\right\}.
\]
\end{lemma}
\begin{proof} This can be proved either by adapting the argument in section 2.4
of \cite{HR1} or by a straightforward computation which exploits the fact that
$\forall h\in SU(2)$, $\overline{h}J = Jh$, $J\epsilon = i\epsilon$ and
$J\overline{\epsilon} = -i \overline{\epsilon}$.  \end{proof}

\noindent We conclude that if, using $\u(3)^{\mathbb{C}} =
\u(3)^{\mathbb{C}}_{-1}\oplus \u(3)^{\mathbb{C}}_0\oplus
\u(3)^{\mathbb{C}}_1\oplus \u(3)^{\mathbb{C}}_2$, we decompose $A$ as
\[
A = A_{-1} + A_0 + A_1 + A_2,
\]
where each $A_a$ is a 1-form with coefficients in the $\u(3)^{\mathbb{C}}_a$,
then we recover the previous splitting by setting $A_0 = A_{\su(2)}$, $A_2 =
A_{\u(1)}$, $A_{-1} = A_{\mathbb{C}^2}'$ and $A_1 = A_{\mathbb{C}^2}''$.  Note
that the two last conditions actually reflects the conformality of $u$.

\begin{remark}
Note that by the automorphism property $[\u(3)_a,\u(3)_b]\subset \u(3)_{a+b\bmod
4}$.
\end{remark}

\subsection{Legendrian framings of Hamiltonian stationary Lagrangian immersions}
\noindent Given the Legendrian framing $F$ of a conformal Lagrangian immersion
$u$ in $\mathbb{C}P^2$, we define the family of deformations $A_\lambda$ of its
Maurer--Cartan form $A$, for $\lambda\in S^1\subset \mathbb{C}^*$ by
\begin{equation}\label{1.6.Alambda}
A_\lambda:= \lambda^{-2}A_2' + \lambda^{-1}A_{-1} + A_0 + \lambda A_1 +
\lambda^2A_2'',
\end{equation}
where $A_2':= A_2(\partial /\partial z)dz$ and $A_2'':= A_2(\partial /\partial
\bar{z})d\bar{z}$.  We then have the following:

\begin{theorem}\label{1.6.theo}
Given a conformal Lagrangian immersion $u:\Omega\longrightarrow \mathbb{C}P^2$
and a Legendrian framing $F$ of $u$, the Maurer--Cartan form of $F$ satisfies
\begin{equation}\label{1.6.lag}
A_{-1} = A_{-1}'= A_{\mathbb{C}^2}'\hbox{ and }A_1 = A_1''= A_{\mathbb{C}^2}'.
\end{equation}
Furthermore $u$ is Hamiltonian stationary if and only if, defining $A_\lambda$
as in (\ref{1.6.Alambda}),
\begin{equation}\label{1.6.curvature}
dA_\lambda + A_\lambda\wedge A_\lambda = 0,\quad \forall \lambda\in S^1.
\end{equation}
\end{theorem}
\begin{remark}
For $\lambda = 1$, $A_1 = A$ and equation (\ref{1.6.curvature}) is a consequence
of its definition $A:=F^{-1}\cdot dF$.
\end{remark}
\begin{proof} See \cite{HR1} and \cite{HR3}.\end{proof}
\noindent We remark that all the conditions that have been collected about the
components $A_a$ can be encoded by the following twisting condition on
$A_\lambda$:
\[
\forall \lambda\in S^1,\quad \tau(A_\lambda) = A_{i\lambda}.
\]
Thus we are led to define the following {\em twisted loop algebra}
\[
\Lambda \u(3)_\tau := \{S^1\ni\lambda \longmapsto \xi_\lambda\in \u(3)/\forall
\lambda\in S^1, \tau(\xi_\lambda) = \xi_{i\lambda}\},
\]
and $A_\lambda$ is a 1-form on $\Omega$ with coefficients in $\Lambda
\u(3)_\tau$.\\

\noindent Actually $\Lambda \u(3)_\tau$ is the Lie algebra of the following {\em
twisted loop group}
\[
\Lambda U(3)_\tau := \{S^1\ni\lambda \longmapsto g_\lambda\in U(3)/\forall
\lambda\in S^1, \tau(g_\lambda) = g_{i\lambda}\},
\]
where the Lie algebra automorphism $\tau :\u(3)^{\mathbb{C}}\longrightarrow
\u(3)^{\mathbb{C}}$ has been extended to the Lie group automorphism by
\[
\begin{array}{cccc}
\tau: & U(3)^{\mathbb{C}} & \longrightarrow & U(3)^{\mathbb{C}}\\
 & M & \longmapsto & \left( \begin{array}{cc}-J & 0\\0 & 1\end{array}\right)
 \cdot \,^t\!M^{-1}\cdot \left( \begin{array}{cc}J & 0\\0 & 1\end{array}\right).
\end{array}
\]
Now if we assume that $\Omega$ is simply connected, then relation
(\ref{1.6.curvature}) allows us to integrate $A_\lambda$, i.e.\,to find a map
$F_\lambda:\Omega\longrightarrow U(3)$ for any $\lambda\in S^1$ such that
$dF_\lambda = F_\lambda\cdot A_\lambda$.  Moreover if we choose some base point
$z_0\in \Omega$, then by requiring further that $F_\lambda(z_0) = Id$,
$F_\lambda$ is unique.  A key observation is then that $\tau(A_\lambda) =
A_{i\lambda}$ implies $\tau(F_\lambda) = F_{i\lambda}$, $\forall \lambda\in
S^1$.  Hence, a conformal Lagrangian immersion $u:\Omega\longrightarrow
\mathbb{C}P^2$ is Hamiltonian stationary if and only if any Legendrian lift $F$
of it can be deformed into a map $F_\lambda:\Omega\longrightarrow \Lambda
U(3)_\tau$, such that $F_\lambda^{-1}\cdot dF_\lambda$ has the form
(\ref{1.6.Alambda}).  Summarizing this result with the observations in the
previous section we have:

\begin{theorem}\label{1.6.theobis}
Given a simply connected domain $\Omega\subset \mathbb{C}$ and a base point
$z_0\in \Omega$, the set of Hamiltonian stationary conformal Lagrangian
immersions $u:\Omega\longrightarrow \mathbb{C}P^2$ such that $u(z_0) = [0:0:1]$
is in bijection with the set of maps $F_\lambda:\Omega\longrightarrow \Lambda
U(3)_\tau$, such that $F_\lambda(z_0)=\mathit{Id}$ and the Fourier decomposition
of $A_\lambda:=F_\lambda^{-1}\cdot dF_\lambda$, $A_\lambda = \sum_{k\in
\mathbb{Z}}\widehat{A}_k\lambda^k$ satisfies
\begin{equation}\label{1.6.-infty}
\forall k\in \mathbb{Z},\quad k\leq -3 \Longrightarrow \widehat{A}_k = 0,
\end{equation}
\begin{equation}\label{1.6.-2}
\widehat{A}_{-2}=
a(z)dz\left(\begin{array}{ccc}i&0&0\\0&i&0\\0&0&0\end{array}\right), \quad
\hbox{where }a\in{\cal C}^\infty(\Omega,\mathbb{C}),
\end{equation}
\begin{equation}\label{1.6.-1}
\widehat{A}_{-1}=\widehat{A}_{-1}(\partial /\partial z)dz,\quad \hbox{i.e. }
\widehat{A}_{-1}(\partial /\partial \bar{z}) = 0.
\end{equation}
\end{theorem}
\begin{proof} For any conformal Lagrangian Hamiltonian stationary immersion $u$
the existence of $F_\lambda$ and the properties (\ref{1.6.-infty}),
(\ref{1.6.-2}) and (\ref{1.6.-1}) are immediate consequences of Theorem
\ref{1.6.theo}.  Conversely for any map $F_\lambda$, conditions
(\ref{1.6.-infty}), (\ref{1.6.-2}) and (\ref{1.6.-1}) and the reality condition
$\overline{A_\lambda} = A_\lambda$ imply that $A_\lambda$ must satisfy
(\ref{1.6.Alambda}).  In particular we remark that condition (\ref{1.6.-2}) is a
reformulation of (\ref{A33=0}).  Thus by theorem \ref{1.6.theo} we deduce that
$F_1$ is the Legendrian lift of some Hamiltonian stationary conformal Lagrangian
immersion.\end{proof}
\begin{remark}
>From the analysis of the Maurer--Cartan form of a Legendrian lift we know that
actually the function $a$ in (\ref{1.6.-2}) is ${1\over 2}\partial
\beta/\partial z$, where $\beta$ is the Lagrangian angle function.  In
particular since $u$ is Hamiltonian stationary $\beta$ is harmonic and hence $a$
is holomorphic.
\end{remark}

\subsection{An alternative characterization}
We introduce here another construction using based loop groups for
characterizing Hamiltonian stationary Lagrangian conformal immersions.  Consider
\[
E_\lambda:= F_\lambda\cdot F^{-1}.
\]
We can observe that $E_\lambda$ is a map with values in the {\em based loop
group}
\[
\Omega U(3):= \{S^1\ni\lambda\longmapsto g_\lambda\in U(3)/g_{\lambda=1}=1\},
\]
since $F_{\lambda=1}=F$.  It is easy to check that $\Omega U(3)$ is a loop
group, the Lie algebra of which is
\[
\Omega \u(3):= \{S^1\ni\lambda\longmapsto \xi_\lambda\in
\u(3)/\xi_{\lambda=1}=0\}.
\]
Note that the (formal) Fourier expansion of an element $\xi_\lambda\in \Omega
\u(3)$ can be written $\xi_\lambda = \sum_{k\in\mathbb{Z}\setminus
\{0\}}\widehat{\xi}_k(\lambda^k-1)$.\\

\noindent The Maurer--Cartan form of $E_\lambda$ is
\[
\begin{array}{ccl}
\Gamma_\lambda & := & E_\lambda^{-1}\cdot dE_\lambda \\
& = & F\cdot\left( F_\lambda^{-1}\cdot dF_\lambda - F\cdot dF\right) \cdot
F^{-1} = F\cdot\left(A_\lambda - A\right) \cdot F^{-1}\\
& = & (\lambda^{-2}-1)\Gamma_2' + (\lambda^{-1}-1)\Gamma_{-1} +
(\lambda-1)\Gamma_1 + (\lambda^{2}-1)\Gamma_2'',
 \end{array}
\]
where $\Gamma_2':= F\cdot A_2'\cdot F^{-1}$, $\Gamma_{-1}:= F\cdot A_{-1}\cdot
F^{-1}$, $\Gamma_1:= F\cdot A_1\cdot F^{-1}$ and $\Gamma_2'':= F\cdot A_2''\cdot
F^{-1}$.  We can observe in particular that
\[
\Gamma_2' = ia \pi^\perp dz,\quad \hbox{where }\pi^\perp:=
F\cdot\left(\begin{array}{ccc}1&&\\&1&\\&&0\end{array}\right)\cdot F^{-1}.
\]
Note that $\pi^\perp$ is the Hermitian orthogonal projection in $\C^3$ onto the
plane $\widehat{u}^\perp$ (moreover $\pi^\perp$ is actually independent of the
lift $\widehat{u}$ chosen for $u$).\\

\noindent Lastly we point out the following equivariance property with respect
to the automorphism $\tau_u$ defined\footnote{ We can remark that the definition
of $\tau_u$ is independent from the choice of the Legendrian framing $F$ of $u$,
and depends only on $u$.  This means that for any pair of Legendrian framings
$F$ and $\widehat{F}$ such that $ \widehat{F} = F\cdot G$, where
$G=\left(\begin{array}{cc}g&\\&1\end{array}\right)$ and $g:\Omega\longrightarrow
SU(2)$, we have $\widehat{F}\cdot \tau(\widehat{F}^{-1}\cdot M\cdot
\widehat{F})\cdot \widehat{F}^{-1} = F\cdot \tau (F^{-1}\cdot M\cdot F)\cdot
F^{-1}$.  This can be checked by a computation using the fact that
$\left(\begin{array}{cc}\pm J&\\&1\end{array}\right)\cdot G =\overline{G}\cdot
\left(\begin{array}{cc}\pm J&\\&1\end{array}\right)$.} by
\[
\tau_u(M) = F\cdot \tau (F^{-1}\cdot M\cdot F)\cdot F^{-1}.
\]
We have obviously $\tau_u^4=1$.  Moreover, setting
\begin{eqnarray*}
\gamma_\lambda & := & \lambda^{-2}\Gamma_2' + \lambda^{-1}\Gamma_{-1} +
\lambda\Gamma_1+ \lambda^{2}\Gamma_2'' \\
& = & F\cdot\left( \lambda^{-2}A_2'+\lambda^{-1}A_{-1} + \lambda A_1+
\lambda^{2}A_2''\right)\cdot F^{-1}
\end{eqnarray*}
and $\gamma:= \gamma_{\lambda=1} = F\cdot\left( A_2'+A_{-1} + A_1+
A_2''\right)\cdot F^{-1}$, so that $\Gamma_\lambda = \gamma_\lambda -\gamma$, we
have
\[
\tau_u(\gamma_\lambda) = \gamma_{i\lambda}.
\]

\section{Finite type solutions}
\noindent In \cite{HR3} we showed how Theorem \ref{1.6.theobis} allows us to
adapt the theory of J. Dorfmeister, F. Pedit and H.Y. Wu \cite{DPW}, in order to
build a Weierstrass type representation theory of {\em all} conformal Lagrangian
Hamiltonian stationary immersions, i.e.\,using holomorphic data.  Here we want
to exploit Theorem \ref{1.6.theobis} in order to construct a particular class of
examples of solutions: the {\em finite type} ones.

\subsection{Definitions}
\noindent We invite the Reader to consult \cite{BFPP}, \cite{G} or \cite{H} for
more details.  We first observe that $U(3)_0:=\{g\in U(3)/ \tau(g)=g\}$, the
fixed set of $\tau$, is a subgroup of $U(3)$, the Lie algebra of which is
$\u(3)_0$ (same observation about $U(3)_0^{\C}$).  Actually $U(3)_0$ is
isomorphic to $SU(2)$ so that we make the identifications $U(3)_0\simeq SU(2)$
and $\u(3)_0\simeq \su(2)$.  We will need an Iwasawa decomposition of
$SU(2)^{\mathbb{C}}$ for our purpose: it will be a pair $(SU(2),\mathfrak{B})$
of {\em subgroups} of $SU(2)^{\mathbb{C}}$, such that $\forall g\in
SU(2)^{\mathbb{C}}$, $\exists !  (f,b)\in SU(2)\times \mathfrak{B}$ with
$g=f\cdot b$, a property that we summarize by writing $SU(2)^{\mathbb{C}} =
SU(2)\cdot \mathfrak{B}$.  Moreover $\mathfrak{B}$ is a solvable Borel subgroup.
We can choose for example
\[
\mathfrak{B}:=\left\{\left(\begin{array}{cc}T^1_1 &
0\\T^2_1&T^2_2\end{array}\right)/ T^1_1,T^2_2\in (0,\infty), T^2_1\in
\mathbb{C}, T^1_1T^2_2 = 1\right\}.
\]
We denote by $\mathfrak{b}$ the Lie algebra of $\mathfrak{B}$.  The Iwasawa
decomposition $SU(2)^{\mathbb{C}} = SU(2)\cdot \mathfrak{B}$ immediately implies
the vector space decomposition $\su(2)^{\mathbb{C}} = \su(2)\oplus
\mathfrak{b}$, which leads to the definition of the two projection mappings
$(\cdot )_{\su}:\su(2)^{\mathbb{C}} \longrightarrow \su(2)$ and
$(\cdot)_\mathfrak{b}:\su(2)^{\mathbb{C}} \longrightarrow \mathfrak{b}$ such
that
\[
\forall \xi \in \su(2)^{\mathbb{C}}, \quad \xi =
(\xi)_{\su}+(\xi)_\mathfrak{b}\quad \hbox{with }(\xi)_{\su}\in \su(2)\hbox{ and
}(\xi)_\mathfrak{b} \in \mathfrak{b}.
\]

\noindent Then we define the following twisted loop algebras
\[
\Lambda \u(3)_\tau^{\mathbb{C}} := \{S^1\ni\lambda \longmapsto \xi_\lambda\in
\u(3)^{\mathbb{C}}/ \forall \lambda\in S^1, \tau(\xi_\lambda) =
\xi_{i\lambda}\},
\]
\[
\Lambda_{\mathfrak{b}}^+ \u(3)_\tau^{\mathbb{C}} := \{[\lambda \longmapsto
\xi_\lambda] \in \Lambda \u(3)_\tau^{\mathbb{C}}/\forall k\in \mathbb{Z}, k\leq
-1\Longrightarrow \widehat{\xi}_k = 0\hbox{ and }\widehat{\xi}_0\in
\mathfrak{b}\},
\]
where we use the Fourier decomposition $\xi_\lambda =\sum_{k\in
\mathbb{Z}}\widehat{\xi}_k\lambda^k$.\\

\noindent The decomposition $\su(2)^{\mathbb{C}} = \su(2)\oplus \mathfrak{b}$
can be extended to loop algebras, i.e. to the splitting $\Lambda
\u(3)_\tau^{\mathbb{C}} = \Lambda \u(3)_\tau\oplus \Lambda_{\mathfrak{b}}^+
\u(3)_\tau^{\mathbb{C}}$.  This can be checked by using the Fourier expansion of
an element $\xi_\lambda\in \Lambda \u(3)_\tau^{\mathbb{C}}$:
\[
\sum_{k\in \mathbb{Z}}\widehat{\xi}_k\lambda^k = \left(
\sum_{k<0}\widehat{\xi}_k\lambda^k + (\widehat{\xi}_0)_{\su} -
\sum_{k>0}\left(\widehat{\xi}_{-k}\right)^\dagger\lambda^k \right) + \left(
(\widehat{\xi}_0)_\mathfrak{b} + \sum_{k>0}\left(\widehat{\xi}_k
+\left(\widehat{\xi}_{-k}\right)^\dagger \right)\lambda^k \right).
\]
We will denote the corresponding projection mappings by $(\cdot
)_{\Lambda_{\su}}:\Lambda \u(3)_\tau^{\mathbb{C}}\longrightarrow \Lambda
\u(3)_\tau$ and $(\cdot)_{\Lambda_\mathfrak{b}^+}:\Lambda
\u(3)_\tau^{\mathbb{C}}\longrightarrow \Lambda_{\mathfrak{b}}^+
\u(3)_\tau^{\mathbb{C}}$.\\

\noindent We also introduce the following finite dimensional subspaces of
$\Lambda \u(3)_\tau$: for any $p\in \mathbb{N}$ we let
\[
\Lambda^{2+4p} \u(3)_\tau:= \left\{[\lambda \longmapsto \xi_\lambda]\in \Lambda
\u(3)_\tau/ \xi_\lambda =\sum_{k=-2-4p}^{2+4p}\widehat{\xi}_k\lambda^k\right\}.
\]
We can now define a pair of vector fields $X_1,X_2:\Lambda^{2+4p}
\u(3)_\tau\longrightarrow \Lambda \u(3)_\tau$ by
\begin{equation}\label{2.1.X}
X_1(\xi_\lambda):= [\xi_\lambda,(\lambda^{4p}\xi_\lambda)_{\Lambda_{\su}}],
\quad X_2(\xi_\lambda):=
[\xi_\lambda,(i\lambda^{4p}\xi_\lambda)_{\Lambda_{\su}}].
\end{equation}
Note that $\lambda^{4p}\xi_\lambda$ belongs to $\Lambda \u(3)_\tau^{\C}$, so
that $(\lambda^{4p}\xi_\lambda)_{\Lambda_{\su}}$ is well defined.
\begin{lemma}
Let $p\in \mathbb{N}$ and $X_1$ and $X_2$ defined by (\ref{2.1.X}).  Then
\begin{itemize}
\item $\forall \xi_\lambda\in \Lambda^{2+4p} \u(3)_\tau$, $X_1(\xi_\lambda),
X_2(\xi_\lambda)\in T_{\xi_\lambda}\Lambda^{2+4p} \u(3)_\tau \simeq
\Lambda^{2+4p} \u(3)_\tau$, so that $X_1$ and $X_2$ are tangent vector fields to
$\Lambda^{2+4p} \u(3)_\tau$.  \item $|\xi_\lambda|^2$ is preserved by $X_1$ and
$X_2$.  Hence the flow of these vector fields are defined for all time \item The
Lie bracket of $X_1$ and $X_2$ vanishes:
\begin{equation}\label{2.1.[,]}
[X_1,X_2] = 0.
\end{equation}
\end{itemize}
\end{lemma}
\begin{proof} This result follows by a straightforward adaptation of the
analogous results for harmonic maps in \cite{BFPP} (see e.g.\,\cite{G} and
\cite{H}).  Note that the proof of (\ref{2.1.[,]}) rests upon the crucial
property that $\Lambda \u(3)_\tau$ and $\Lambda_{\mathfrak{b}}^+
\u(3)_\tau^{\mathbb{C}}$ are Lie algebras (see e.g.\,\cite{BP},
\cite{H}).\end{proof}

\noindent This result allows us to integrate simultaneously $X_1$ and $X_2$.  So
for any $\xi_\lambda^0\in \Lambda^{2+4p} \u(3)_\tau$ there exists a unique map
$\xi_\lambda:\mathbb{R}^2\longrightarrow \Lambda^{2+4p} \u(3)_\tau$ such that
$\xi_\lambda(z_0) = \xi_\lambda^0$ and
\begin{equation}\label{2.1.flot}
{\partial \xi_\lambda\over \partial x}(x,y) = X_1\left(\xi_\lambda(x,y)\right)
\quad \hbox{and} \quad {\partial \xi_\lambda\over \partial y}(x,y) =
X_2\left(\xi_\lambda(x,y)\right).
\end{equation}
Denoting by $z=x+iy\in \C$, the system (\ref{2.1.flot}) can be rewritten
\[
\begin{array}{ccl}
d\xi_\lambda & = & \left[\xi_\lambda,
\left(\lambda^{4p}\xi_\lambda\right)_{\Lambda_{\su}}dx +
\left(i\lambda^{4p}\xi_\lambda\right)_{\Lambda_{\su}}dy\right]\\
& = & \left[\xi_\lambda, \left(\lambda^{4p}\xi_\lambda
dz\right)_{\Lambda_{\su}}\right].
\end{array}
\]
Let us denote by $A_\lambda:= (\lambda^{4p}\xi_\lambda dz)_{\Lambda_{\su}}$.
Since the system (\ref{2.1.flot}) is overdetermined, $A_\lambda$ should satisfy
a compatibility condition.  Indeed one can check that
\begin{equation}\label{2.1.curv}
dA_\lambda + A_\lambda\wedge A_\lambda = 0.
\end{equation}
This relation can be proved by a method similar to the proof of (\ref{2.1.[,]})
(see \cite{H}).  It implies that there exists a map $F_\lambda:\C\longrightarrow
\Lambda U(3)_\tau$ such that
\begin{equation}\label{2.1.dF}
dF_\lambda = F_\lambda\cdot A_\lambda.
\end{equation}
Now observe that $\lambda^{4p}\xi_\lambda =
\sum_{k=-2}^{8p+2}\hat{\xi}_{k-4p}\lambda^k$ implies
\[
A_\lambda = \lambda^{-2}\hat{\xi}_{-4p-2}dz + \lambda^{-1}\hat{\xi}_{-4p-1}dz +
\left(\hat{\xi}_{-4p}dz\right)_{\su} -
\lambda\left(\hat{\xi}_{-4p-1}\right)^\dagger d\bar{z} -
\lambda^{2}\left(\hat{\xi}_{-4p-2}\right)^\dagger d\bar{z}.
\]
We recall that $\hat{\xi}_{-4p-2}\in \u(3)_2^{\C}$ and so has the form
$\mathrm{diag}(ia,ia,ib)$.  Moreover we have the following result.
\begin{lemma}
If $\xi_\lambda\longrightarrow \Lambda^{2+4p} \u(3)_\tau$ and $A_\lambda:=
(\lambda^{4p}\xi_\lambda dz)_{\Lambda_{\su}}$ are solutions of $d\xi_\lambda =
[\xi_\lambda,A_\lambda]$, then $\hat{\xi}_{-4p-2}$ is constant.
\end{lemma}
\begin{proof} The relevant term in the Fourier expansion of $d\xi_\lambda =
[\xi_\lambda,A_\lambda]$ gives
\[
\begin{array}{ccl}
\displaystyle d\hat{\xi}_{-4p-2} & = & \displaystyle \left[ \hat{\xi}_{-4p-2},
\left(\hat{\xi}_{-4p}dz\right)_{\su}\right] + \left[ \hat{\xi}_{-4p-1},
\hat{\xi}_{-4p-1}\right]dz + \left[ \hat{\xi}_{-4p},
\hat{\xi}_{-4p-2}\right]dz\\
 & = & \displaystyle
 \left[\left(\hat{\xi}_{-4p}dz\right)_\mathfrak{b},\hat{\xi}_{-4p-2}\right].
\end{array}
\]
But since the coefficients of $\left(\hat{\xi}_{-4p}dz\right)_\mathfrak{b}$ are
in $\u(3)_0^{\C}$ and $\hat{\xi}_{-4p-2}$ takes values in $\u(3)_2$ we deduce
that $d\hat{\xi}_{-4p-2} = 0$, because $\u(3)_0^{\C}$ and $\u(3)_2^{\C}$
commute.  \end{proof}

We deduce from this result that if we choose the initial value $\xi_\lambda^0$
of $\xi_\lambda$ to be such that $\hat{\xi}^0_{-4p-2} = \mathrm{diag}(ia,ia,0)$
then $\hat{\xi}_{-4p-2}$ is equal to that value for all $(x,y)$.  So in this
case the map $F_\lambda$ obtained by integrating $A_\lambda$ satisfies all the
requirements of Theorem \ref{1.6.theobis}.  It implies that $F_\lambda$
represents a (conjugate family) of Hamiltonian stationary conformal Lagrangian
immersion(s).  The category of such $F_\lambda$'s are exactly characterized by
the following definition.

\begin{definition}\label{2.1.def-ft}
Let $F_\lambda$ be a family of Hamiltonian stationary conformal Lagrangian
immersions and let $A_\lambda:= F_\lambda^{-1}\cdot dF_\lambda$.  Then
$F_\lambda$ is called a family of {\em finite type solutions} if and only if
there exists $p\in \mathbb{N}$ and a map $\xi_\lambda:\mathbb{C}\longrightarrow
\Lambda^{2+4p} \u(3)_\tau$ such that $\hat{\xi}_{-4p-2} =
\mathrm{diag}(ia,ia,0)$, for some constant $a\in \C$, and
\begin{equation}\label{2.1.laxfn}
d\xi_\lambda = [\xi_\lambda,A_\lambda]
\end{equation}
\begin{equation}\label{2.1.laxfn2}
(\lambda^{4p}\xi_\lambda dz)_{\Lambda_{\su}} = A_\lambda.
\end{equation}
\end{definition}
We also need the following definition in which we introduce an a priori weaker
notion of finite type solution.

\begin{definition}
Let $F_\lambda$ be a family of Hamiltonian stationary conformal Lagrangian
immersions and let $A_\lambda:= F_\lambda^{-1}\cdot dF_\lambda$.  Then
$F_\lambda$ is called a family of {\em quasi-finite type solutions} if and only
if it satisfies the same requirements as in definition \ref{2.1.def-ft} excepted
that condition (\ref{2.1.laxfn2}) is replaced by
\begin{equation}\label{2.1.laxqfn2}
\exists B\in \Omega^1\otimes \u(3)^{\C}_0,\quad (\lambda^{4p}\xi_\lambda
dz)_{\Lambda_{\su}} = A_\lambda + B.
\end{equation}
\end{definition}
We shall see in Section 3.3 that both definitions are actually equivalent.

\subsection{An alternative description of quasi-finite type solutions}
We may as well characterize such finite type solutions in terms of $E_\lambda =
F_\lambda\cdot F^{-1}$.  For that purpose we need to introduce the untwisted
loop Lie algebra
\[
\Lambda^+\u(3)^{\C}:= \{S^1\ni\lambda\longmapsto \xi_\lambda\in \u(3)^{\C}/
\xi_\lambda = \sum_{k=0}^\infty \widehat{\xi}_k\lambda^k\}
\]
and observe that any $\xi_\lambda=\sum_{k=-\infty}^\infty
\widehat{\xi}_k\lambda^k\in \Lambda \u(3)^{\C}$ can be split as
\[
\xi_\lambda = \left(\sum_{k=-\infty}^{-1}\widehat{\xi}_k(\lambda^k-1) -
(\widehat{\xi}_k)^\dagger (\lambda^{-k}-1)\right) + \left(\sum_{k=0}^\infty
\widehat{\xi}_k\lambda^k + \sum_{k=1}^\infty \widehat{\xi}_{-k} +
(\widehat{\xi}_{-k})^\dagger (\lambda^{k}-1)\right)
\]
and hence $\Lambda \u(3)^{\C} = \Omega \u(3)\oplus \Lambda^+\u(3)^{\C}$.  This
defines a pair of projection mappings $(\cdot)_{\Omega}:\Lambda \u(3)^{\C}
\longrightarrow \Omega \u(3)$ and $(\cdot)_{\Lambda^+}:\Lambda \u(3)^{\C}
\longrightarrow \Lambda^+\u(3)^{\C}$.

Now consider a family $F_\lambda$ of quasi-finite type, let $A_\lambda:=
F_\lambda^{-1}\cdot dF_\lambda$, $A:= F^{-1}\cdot dF$ (where $F=F_{\lambda=1}$)
and $\xi_\lambda$ be a solution of (\ref{2.1.laxfn}).  We let
\[
\eta_\lambda:= F\cdot \xi_\lambda\cdot F^{-1} = \sum_{k=-2-4p}^{2+4p}F\cdot
\widehat{\xi}_k\cdot F^{-1}\lambda^k .
\]
Then (\ref{2.1.laxfn}) implies by a straightforward computation that
\[
\begin{array}{ccl}
d\eta_\lambda & = & F\cdot (d\xi_\lambda + [A,\xi_\lambda])\cdot F^{-1}\\
& = & F\cdot ([\xi_\lambda,A_\lambda] - [\xi_\lambda,A])\cdot F^{-1} =
[\eta_\lambda, \Gamma_\lambda],
\end{array}
\]
where $\Gamma_\lambda:= E_\lambda^{-1}\cdot dE_\lambda$.  Now setting
$R_\lambda:= \sum_{k=-4p}^{2+4p}F\cdot \widehat{\xi}_k\cdot F^{-1}\lambda^k$, we
have
\begin{eqnarray*}
\displaystyle \left(\lambda^{4p}\eta_\lambda dz\right)_\Omega & = &
\displaystyle \left( \lambda^{-2}F\cdot \widehat{\xi}_{-2-4p}\cdot F^{-1}dz +
\lambda^{-1}F\cdot \widehat{\xi}_{-1-4p}\cdot F^{-1}dz + \lambda^{4p}R_\lambda
dz\right)_\Omega\\
& = & \displaystyle (\lambda^{-2}-1)F\cdot \widehat{\xi}_{-2-4p}\cdot F^{-1}dz +
(\lambda^{-1}-1)F\cdot \widehat{\xi}_{-1-4p}\cdot F^{-1}dz\\
&& \displaystyle - (\lambda-1)\left(F\cdot \widehat{\xi}_{-1-4p}\cdot
F^{-1}\right)^\dagger d\bar{z} - (\lambda^{2}-1)\left(F\cdot
\widehat{\xi}_{-2-4p}\cdot F^{-1}\right)^\dagger d\bar{z}.
\end{eqnarray*}
But relation (\ref{2.1.laxqfn2}) implies in particular that
$\widehat{\xi}_{-2-4p} = A_2'(\partial /\partial z)$ and $\widehat{\xi}_{-1-4p}
= A_{-1}(\partial /\partial z)$.  So we deduce that
\[
\displaystyle \left(\lambda^{4p}\eta_\lambda dz\right)_\Omega = F\cdot
\left(A_\lambda - A\right) \cdot F^{-1} = \Gamma_\lambda.
\]
Hence $E_\lambda$ can be constructed by solving a system analogous to
(\ref{2.1.laxfn}), (\ref{2.1.laxfn2}), i.e.
\begin{equation}\label{2.2.laxfn}
d\eta_\lambda + \left[\Gamma_\lambda,\eta_\lambda\right] = 0 \quad \hbox{and}
\quad \Gamma_\lambda = \left(\lambda^{4p}\eta_\lambda dz\right)_\Omega.
\end{equation}
Conversely a similar computation shows that a solution of (\ref{2.2.laxfn})
gives rise to a quasi-finite type family of solutions by an inverse
transformation, but we shall prove more in the next section.\\

\noindent Note that system (\ref{2.2.laxfn}) can also be interpreted as a pair
of commuting ordinary differential equations in the finite dimensional space
$\Lambda^{2+4p}\u(3):= \{S^1\ni\lambda\longmapsto\eta_\lambda \in \u(3)/
\eta_\lambda =\sum_{k=-2-4p}^{2+4p}\widehat{\eta}_k\lambda^k\}$.  It is the
analogue of the definition of a finite type solution according to \cite{BFPP}.\\

\subsection{Quasi-finite type solutions are actually finite type}

\noindent We show here the following
\begin{theorem}\label{3.3.qtf=tf}
For any family $F_\lambda$ of Hamiltonian stationary Lagrangian conformal
immersions of quasi-finite type, i.e.\,such that there exists
$\xi_\lambda:\Omega\longrightarrow \Lambda^{2+4p}\u(3)_\tau$ which satisfies
(\ref{2.1.laxfn}) and (\ref{2.1.laxqfn2}), there exists a gauge transformation
$F_\lambda\longmapsto F_\lambda^G:=F_\lambda\cdot G$, where $G\in {\cal
C}^\infty(\Omega,U(3)_0)$, such that $F_\lambda^G$ is of finite type.  More
precisely, denoting by $A_\lambda^G:=G^{-1}\cdot A_\lambda\cdot G + G^{-1}\cdot
dG$ and $\xi_\lambda^G:=G^{-1}\cdot \xi_\lambda\cdot G$, then $d\xi_\lambda^G
+\left[A_\lambda^G,\xi_\lambda^G\right] = G^{-1}\cdot \left(d\xi_\lambda
+[A_\lambda,\xi_\lambda]\right)\cdot G = 0$ and $(\lambda^{4p}\xi_\lambda^G
dz)_{\Lambda_{\su}} =A_\lambda^G$.
\end{theorem}
\begin{proof} We set $E_\lambda:=F_\lambda\cdot F^{-1}$,
$\Gamma_\lambda:=E_\lambda^{-1}\cdot dE_\lambda$ and $\eta_\lambda:= F\cdot
\xi_\lambda\cdot F^{-1}$ and will use the results of the previous section.\\

{\em A constant in $\Lambda^{2+4p}\u(3)_\tau$ associated to the quasi-finite
type family} --- First (\ref{2.2.laxfn}), which is a reformulation of
(\ref{2.1.laxfn}), implies
\[
d\left(E_\lambda\cdot \eta_\lambda\cdot E_\lambda^{-1}\right) = E_\lambda\cdot
\left( d\eta_\lambda +[\Gamma_\lambda,\eta_\lambda]\right)\cdot E_\lambda^{-1}
=0.
\]
Hence
\[
\eta_\lambda^0:= E_\lambda\cdot \eta_\lambda\cdot E_\lambda^{-1}
\]
is a constant in $\Lambda \u(3)$.  Moreover
\[
\eta_\lambda^0 = E_\lambda(z_0)\cdot \eta_\lambda(z_0)\cdot E_\lambda^{-1}(z_0)
= \eta_\lambda(z_0) = F(z_0)\cdot \xi_\lambda(z_0)\cdot F^{-1}(z_0)
=\xi_\lambda(z_0),
\]
which proves that $\eta_\lambda^0\in \Lambda^{2+4p}\u(3)_\tau$.\\

{\em An auxiliary map into $\Lambda^+U(3)^{\C}$} --- We let
\[
\Theta_\lambda:= \left(\lambda^{4p}\eta_\lambda dz\right)_{\Lambda^+} =
\lambda^{4p}\eta_\lambda dz - \left(\lambda^{4p}\eta_\lambda dz\right)_\Omega.
\]
Then using (\ref{2.2.laxfn}) we have $\Theta_\lambda = \lambda^{4p}\eta_\lambda
dz - \Gamma_\lambda$ and so
\[
d\Gamma_\lambda + \Gamma_\lambda\wedge \Gamma_\lambda + d\Theta_\lambda -
\Theta_\lambda\wedge \Theta_\lambda = -\lambda^{4p} \left(d\eta_\lambda
+[\Gamma_\lambda,\eta_\lambda]\right) \left({\partial \over \partial
\bar{z}}\right)dz\wedge d\bar{z} =0.
\]
But since $d\Gamma_\lambda + \Gamma_\lambda\wedge \Gamma_\lambda = 0$ this
implies that $d\Theta_\lambda - \Theta_\lambda\wedge \Theta_\lambda = 0$.  Hence
$\exists!  V_\lambda:\Omega\longrightarrow \Lambda^+U(3)^{\C}$ such that
\[
dV_\lambda = \Theta_\lambda\cdot V_\lambda \quad \hbox{and}\quad V_\lambda (z_0)
= 1.
\]
Now, starting from $\lambda^{4p}\eta_\lambda dz =
\Gamma_\lambda+\Theta_\lambda$, we deduce that
\[
\begin{array}{ccl}
\displaystyle \lambda^{4p}\eta_\lambda^0dz & = & \displaystyle E_\lambda\cdot
\Gamma_\lambda\cdot E_\lambda^{-1} + E_\lambda\cdot \Theta_\lambda\cdot
E_\lambda^{-1}\\
& = & \displaystyle dE_\lambda\cdot E_\lambda^{-1} + E_\lambda\cdot
dV_\lambda\cdot V_\lambda^{-1}\cdot E_\lambda^{-1}\\
& = & \displaystyle d\left(E_\lambda\cdot V_\lambda\right)\left(E_\lambda\cdot
V_\lambda\right)^{-1},
\end{array}
\]
which can be integrated into the relation
\[
E_\lambda\cdot V_\lambda = e^{\lambda^{4p}(z-z_0)\eta_\lambda^0}.
\]

{\em An Iwasawa decomposition of $e^{\lambda^{4p}(z-z_0)\eta_\lambda^0}$} ---
The latter implies
\[
e^{\lambda^{4p}(z-z_0)\eta_\lambda^0} = F_\lambda\cdot F^{-1}\cdot V_\lambda.
\]
>From this relation and the fact that $\eta_\lambda^0$ and $F_\lambda$ are
twisted we deduce that $W_\lambda:= F^{-1}\cdot V_\lambda$ is twisted.  It is
also a map with values in $\Lambda^+U(3)^{\C}_\tau$.  However it may not be not
in $\Lambda^+_{\mathfrak{B}}U(3)^{\C}_\tau$ in general, because in the
development
\[
F^{-1}\cdot V_\lambda = \widehat{W}_0 + \sum_{k=1}^\infty
\widehat{W}_k\lambda^k,
\]
we are not sure that $\widehat{W}_0$ takes values in $\mathfrak{B}$.  But it
takes values in $U(3)^{\C}_0$, so by using the Iwasawa decomposition
$U(3)^{\C}_0 = U(3)_0\cdot \mathfrak{B}$ we know that $\exists !  G\in U(3)_0$,
$\exists !\widehat{B}_0\in \mathfrak{B}$, $\widehat{W}_0 = G\cdot
\widehat{B}_0$.  Hence
\[
G^{-1}\cdot F^{-1}\cdot V_\lambda = \widehat{B}_0 + \sum_{k=1}^\infty
G^{-1}\cdot \widehat{W}_k\lambda^k
\]
takes values in $\Lambda^+_{\mathfrak{B}}U(3)^{\C}_\tau$.  So the splitting
\[
e^{\lambda^{4p}(z-z_0)\eta_\lambda^0} = \left(F_\lambda\cdot G\right)
\left(G^{-1}\cdot F^{-1}\cdot V_\lambda\right)
\]
exactly reproduces the Iwasawa decomposition $\Lambda U(3)_\tau^{\C} = \Lambda
U(3)_\tau \cdot \Lambda^+_{\mathfrak{B}}U(3)^{\C}_\tau$ proved in \cite{DPW}.\\

\noindent {\em Conclusion} --- Let us denote by $F_\lambda^G:= F_\lambda\cdot
G$, $A_\lambda^G:= \left(F_\lambda^G\right)^{-1}\cdot dF_\lambda^G = G^{-1}\cdot
A_\lambda\cdot G + G^{-1}\cdot dG$ and $B_\lambda^G:= G^{-1}\cdot F^{-1}\cdot
V_\lambda$ and let us introduce
\[
\xi_\lambda^G:= \left(F_\lambda^G\right)^{-1}\cdot \eta_\lambda^0\cdot
F_\lambda^G.
\]
(These definitions imply immediately $d\xi_\lambda^G +
[A_\lambda^G,\xi_\lambda^G] = 0$.)  The first main observation is that the
relation $\eta_\lambda^0 = E_\lambda\cdot \eta_\lambda\cdot E_\lambda^{-1} =
F_\lambda\cdot F^{-1}\cdot \eta_\lambda\cdot F\cdot F_\lambda^{-1} =
F_\lambda\cdot \xi_\lambda\cdot F_\lambda^{-1}$ implies
\begin{equation}\label{2.3.tf1}
\xi_\lambda^G = G^{-1}\cdot F_\lambda^{-1}\cdot \eta_\lambda^0\cdot
F_\lambda\cdot G = G^{-1}\cdot \xi_\lambda\cdot G.
\end{equation}
Second, from the relation
\[
\begin{array}{ccl}
\displaystyle \lambda^{4p}\eta^0_\lambda dz & = & \displaystyle
d\left(e^{\lambda^{4p}(z-z_0)\eta_\lambda^0}\right) \cdot
e^{-\lambda^{4p}(z-z_0)\eta_\lambda^0}\\
 & = & \displaystyle d\left(F_\lambda^G\cdot B_\lambda^G\right)\cdot
 \left(F_\lambda^G\cdot B_\lambda^G\right)^{-1}\\
& = & \displaystyle dF_\lambda^G\cdot \left(F_\lambda^G\right)^{-1} +
F_\lambda^G\cdot dB_\lambda^G\cdot \left(B_\lambda^G\right)^{-1}\cdot
\left(F_\lambda^G\right)^{-1},
\end{array}
\]
we deduce that
\[
\lambda^{4p}\xi_\lambda^G dz = \left(F_\lambda^G\right)^{-1} \cdot \left(
\lambda^{4p}\eta_\lambda^0 dz\right)\cdot F_\lambda^G =
\left(F_\lambda^G\right)^{-1}\cdot dF_\lambda^G + dB_\lambda^G\cdot
\left(B_\lambda^G\right)^{-1}.
\]
Hence, since $B_\lambda^G$ takes values in
$\Lambda^+_{\mathfrak{B}}U(3)^{\C}_\tau$,
\begin{equation}\label{2.3.tf2}
A_\lambda^G = \left(F_\lambda^G\right)^{-1}\cdot dF_\lambda^G = \left(
\lambda^{4p}\xi_\lambda^Gdz\right)_{\Lambda_{\su}}.
\end{equation}
And relations (\ref{2.3.tf1}) and (\ref{2.3.tf2}) lead to the conclusion.
\end{proof}

\section{All Hamiltonian stationary Lagrangian tori are of finite type}
The subject of this section is to prove the following:
\begin{theorem}\label{4.MainTheo}
Let $u:\mathbb{C}\longrightarrow \mathbb{C}P^2$ be a doubly periodic Hamiltonian
stationary Lagrangian conformal immersion.  Then $u$ is of finite type.
\end{theorem}
We will actually prove a slightly more general result, since we can replace the
doubly periodicity assumption by the hypothesis that the Maurer--Cartan form of
any {\em Legendrian framing} of $u$ is doubly periodic.  This result of course
implies immediately that Hamiltonian stationary Lagrangian tori are of finite
type, since they always can be covered conformally by the plane.\\

\noindent Note also that the study of Hamiltonian stationary Lagrangian tori
splits into exactly two subcases: the {\em minimal} Lagrangian tori and the {\em
non minimal} Hamiltonian stationary Lagrangian ones.  The first case occurs when
the Lagrangian angle function along any Legendrian lift is locally constant, the
second one when this function is harmonic and non constant.  In the case of
minimal Lagrangian surfaces, Theorem \ref{4.MainTheo} is a special case of the
result in \cite{BFPP}, since in this case $u$ is a harmonic map into
$\mathbb{C}P^2$, as discussed in \cite{McI1}, \cite{McI2}, \cite{McI3} and
\cite{J}.  The non minimal case however is not covered by the theory in
\cite{BFPP} and is the subject of this section.  \\

\noindent Let $F:\mathbb{C}\longrightarrow U(3)$ be a Legendrian framing of $u$,
$A:=F^{-1}\cdot dF$ its Maurer--Cartan form and $A_\lambda$ the family of
deformations of $A$ as defined by (\ref{1.6.Alambda}).  The first basic
observation is that $A_2\left({\partial \over \partial z}\right)$ is holomorphic
and doubly periodic on $\C$, hence constant.  Thus two cases occur: either
$A_2\left({\partial \over \partial z}\right)=0$, which corresponds to the
minimal case that we exclude here, or $A_2\left({\partial \over \partial
z}\right)$ is a constant different from 0, the case that we consider next.\\

\noindent In order to show Theorem \ref{4.MainTheo} we need to prove that there
exists some $p\in \mathbb{N}$ and a map $\xi_\lambda:\mathbb{C}\longrightarrow
\Lambda^{2+4p}\u(3)_\tau$ such that $d\xi_\lambda = [\xi_\lambda,A_\lambda]$ and
$A_\lambda = (\lambda^{4p}\xi_\lambda dz)_{\Lambda_{\su}}$.  But thanks to
Theorem \ref{3.3.qtf=tf} it will enough to prove that $A_\lambda -
(\lambda^{4p}\xi_\lambda dz)_{\Lambda_{\su}}$ is a 1-form with coefficients in
$\u(3)_0$.  Our proof here follows a strategy inspired from \cite{BFPP}: a first
step consists in building a formal series $Y_\lambda = \sum_{k=-2}^\infty
\widehat{Y}_k\lambda^k$ which is a solution of $dY_\lambda =
[Y_\lambda,A_\lambda]$.  Such a series is called a {\em formal Killing field}.
We will also require $Y_\lambda$ to be {\em quasi-adapted}, i.e.\,is such that
\begin{equation}\label{4.adapted}
(Y_\lambda dz)_{\Lambda_{\su}} = A_\lambda + B,\quad \hbox{where the
coefficients of }B\hbox{ are in } \u(3)_0.
\end{equation}
This is achieved through a recursion procedure.\\

\noindent In a second step we will show that the coefficients of $Y_\lambda$
form a countable collection of doubly periodic functions satisfying an elliptic
PDE and hence, by using a compactness argument, we conclude that they are
contained in a finite dimensional space.  Then we deduce the existence of
$\xi_\lambda$ using linear algebra.

\subsection{Construction of an adapted formal Killing field}

We first introduce some notations.  We denote by
\[
\pi_0^\perp := \left(\begin{array}{ccc} 1 &&\\ & 1& \\ && 0 \end{array}\right)
\]
and $a := {1\over 2}{\partial \beta\over \partial z}$ (here a constant different
from 0).  Then $A_2' = ia\pi_0^\perp dz$.  We will also set
$X:=A_{-1}\left({\partial \over \partial z}\right)$ and $C:=A_0\left({\partial
\over \partial z}\right)$, so that
\[
A_\lambda = \lambda^{-2}ia\pi_0^\perp dz + \lambda^{-1}Xdz + Cdz -C^\dagger
d\bar{z} -\lambda X^\dagger d\bar{z} + \lambda^2i\overline{a}\pi_0^\perp
d\bar{z}.
\]
We also introduce the linear map $\ad \pi_0^\perp : \u(3)^{\C}\longrightarrow
\u(3)^{\C}$, acting by $\xi\longmapsto [\pi_0^\perp ,\xi]$\footnote{Actually the
map $\xi\longmapsto i[\pi_0^\perp ,\xi]$ corresponds to the complex structure on
the Legendrian distribution.}.  We observe that $\pi_0^\perp $ commutes with the
elements in $\u(3)_0^{\C}$ and $\u(3)_2^{\C}$.  Moreover
\[
\forall a,b\in \C, \; \left[\pi_0^\perp ,\left(\begin{array}{ccc}&& a \\ && b \\
\mp ib &\pm ia \end{array}\right)\right] = \left(\begin{array}{ccc} && a \\ && b
\\ \pm ib & \mp ia& \end{array}\right),
\]
that is $\ad \pi_0^\perp $ maps $\u(3)_{\mp 1}^{\C}$ to $\u(3)_{\pm 1}^{\C}$.
>From that we deduce that $V:= \Ker \ad \pi_0^\perp $ coincides with
$\u(3)_0^{\C}\oplus \u(3)_2^{\C}$ and $V^\perp:= \Im \ad \pi_0^\perp $ coincides
with $\u(3)_{-1}^{\C}\oplus \u(3)_1^{\C}$ (note that $V^\perp$ is actually the
orthogonal subspace to $V$ in $\u(3)^{\C}$).  In our construction we will use
extensively the following properties:
\begin{itemize}
\item the map $\ad \pi_0^\perp \at{V^\perp \to V^\perp}$ is a vector space
isomorphism (it is actually a involution on $V^\perp$), \item the inclusions $VV
\subset V$, $VV^\perp \subset V^\perp$, $V^\perp V \subset V^\perp$ and $V^\perp
V^\perp \subset V$.  These properties can be checked by a direct computation
using the fact that matrices in $V$ are diagonal by blocks and the matrices in
$V^\perp$ are off-diagonal by blocks.  (The three first properties can also be
deduced from the definition of $V$ and $V^\perp$ and the fact that $\ad$ is a
derivation).
\end{itemize}
We look for a formal Killing field $Y_\lambda$, i.e.\,a solution of the equation
\begin{equation} \label{eq:champ_formel}
dY_\lambda = [Y_\lambda,A_\lambda],
\end{equation}
of the form $Y_\lambda = (1+W_\lambda)^{-1}\lambda^{-2}ia\pi_0^\perp
(1+W_\lambda)$, where $W_\lambda = \sum_{k=0}^\infty \widehat{W}_k\lambda^k$ as
in \cite{BFPP}.  In order to have a well-posed problem (and in particular to
guarantee the existence of an unique solution of this type) we assume that
$W_\lambda$ takes values in $V^\perp$.  We start by evaluating
(\ref{eq:champ_formel}) along $\partial/\partial z$.  It gives, after
conjugation by $1+W_\lambda$:
\begin{multline} \label{eq:W}
\lambda^{-2}\frac{\partial a}{\partial z} \pi_0^\perp + \lambda^{-2}
a\big[\pi_0^\perp , \\
 \frac{\partial W_\lambda}{\partial z}(1+W_\lambda)^{-1}
 -(1+W_\lambda)(\lambda^{-2}ia\pi_0^\perp
 +\lambda^{-1}X+C)(1+W_\lambda)^{-1}\big] =0
\end{multline}
Here the fact that $a$ is a constant leads to an immediate simplification,
namely that the bracket in the left hand side of (\ref{eq:W}) is 0.  Thus
equation (\ref{eq:W}) implies that $\frac{\partial W_\lambda}{\partial z}
(1+W_\lambda)^{-1}-(1+W_\lambda)(\lambda^{-2}ia\pi_0^\perp +\lambda^{-1}X+C)
(1+W_\lambda)^{-1}$ lies in $V$, hence there exists a map
$\varphi_\lambda:\C\longrightarrow V$ such that
\[
\pa{W_\lambda}{z}(1+W_\lambda)^{-1} -(1+W_\lambda)(\lambda^{-2}ia\pi_0^\perp
+\lambda^{-1}X+C)(1+W_\lambda)^{-1} = \varphi_\lambda
\]
or
\[
\pa{W_\lambda}{z}-(1+W_\lambda)(\lambda^{-2}ia\pi_0^\perp +\lambda^{-1}X+C) =
\varphi_\lambda(1+W_\lambda),
\]
which can be projected according to the splitting $V \oplus V^\perp$ as
\[
\left\{ \begin{array}{l}\displaystyle \lambda^{-2}ia\pi_0^\perp
+\lambda^{-1}W_\lambda X+C = -\varphi_\lambda \in V
\\
\displaystyle \pa{W_\lambda}{z} - \lambda^{-2}ia W_\lambda \pi_0^\perp -
\lambda^{-1}X - W_\lambda C = \varphi_\lambda W_\lambda \in V^\perp.
\end{array} \right.
\]
Substituting $\varphi_\lambda$,
\[
\pa{W_\lambda}{z} - \lambda^{-2}ia W_\lambda \pi_0^\perp - \lambda^{-1}X-
W_\lambda C + \lambda^{-2}ia\pi_0^\perp W_\lambda+\lambda^{-1}W_\lambda X
W_\lambda +C W_\lambda = 0
\]
or
\[
[ia\pi_0^\perp ,W_\lambda] + \lambda(W_\lambda X W_\lambda-X) +
\lambda^2[C,W_\lambda] + \lambda^2 \pa{W_\lambda}{z} = 0
\]
or
\begin{multline*}
ia\sum_{n\geq 0} [\pi_0^\perp ,\widehat{W}_n]\lambda^n + \sum_{n\geq
1}\left(\sum_{k=0}^{n-1} \widehat{W}_k X \widehat{W}_{n-1-k}\right) \lambda^n -
\lambda X \\
+ \sum_{n\geq 2} \left([C,\widehat{W}_{n-2}] + \pa{\widehat{W}_{n-2}}{z} \right)
\lambda^n = 0.
\end{multline*}
Hence
\[
\left\{ \begin{array}{ll} \displaystyle n=0, & \displaystyle ia[\pi_0^\perp
,\widehat{W}_0] = 0
\\
\displaystyle n=1, & \displaystyle ia[\pi_0^\perp ,\widehat{W}_1]+\widehat{W}_0
X \widehat{W}_0 - X = 0
\\
\displaystyle n \geq 2, & \displaystyle ia[\pi_0^\perp
,\widehat{W}_n]+\sum_{k=0}^{n-1} \widehat{W}_k X \widehat{W}_{n-1-k} +
[C,\widehat{W}_{n-2}] + \pa{\widehat{W}_{n-2}}{z} =0 \end{array} \right.
\]
and thus
\[
\left\{ \begin{array}{l} \widehat{W}_0 = 0
\\
\widehat{W}_1 = -ia^{-1}[\pi_0^\perp ,X]
\\
\displaystyle \widehat{W}_n = ia^{-1}\left[\pi_0^\perp ,\sum_{k=0}^{n-1}
\widehat{W}_k X \widehat{W}_{n-1-k} + [C,\widehat{W}_{n-2}] +
\pa{\widehat{W}_{n-2}}{z}\right] \end{array} \right.
\]
We observe that the formal Killing field is {\em quasi-adapted} in the sense
that the two first coefficients are the right ones:
\begin{eqnarray*}
Y_\lambda &=& ia\lambda^{-2} (1+W_\lambda)^{-1}\pi_0^\perp (1+W_\lambda) = ia
\lambda^{-2} \left( \pi_0^\perp - \lambda [\widehat{W}_1,\pi_0^\perp ] +
\O(\lambda^2) \right) \\
&=& \lambda^{-2} i a \pi_0^\perp + \lambda^{-1} X + \O(1).
\end{eqnarray*}
Another pleasant property is that this formal field is automatically twisted (as
in the case of $\C^2$, see \cite{HR1}).  Indeed using the fact that $\tau$ is an
automorphism for the product of matrices as well as for the Lie bracket (and so
$[\u(3)^{\C}_a,\u(3)^{\C}_b] \subset \u(3)^{\C}_{a+b}$ and $\u(3)^{\C}_a
\u(3)^{\C}_b \subset \u(3)^{\C}_{a+b}$), we obtain that
\[
\tau(Y_\lambda) = \tau(1+W_\lambda)^{-1} (-\lambda^2)ia\pi_0^\perp
\tau(1+W_\lambda).
\]
Thus it is enough to show that $1+W_\lambda$ is twisted, i.e.\,$\tau(W_\lambda)
= W_{i\lambda}$.  In terms of the Fourier decomposition of $W_\lambda$ this is
equivalent to proving that $\widehat{W}_n$ belongs to $\u(3)^{\C}_n$.  Let us
prove it by recursion.  We already know that $\widehat{W}_1 =
-ia^{-1}[\pi_0^\perp ,X]$ is in $\u(3)^{\C}_1$.  Assume that the result is true
up to $n-1$,then
\[
\sum_{p=0}^{n-1} \widehat{W}_p X \widehat{W}_{n-1-p} + [C,\widehat{W}_{n-2}] +
\pa{\widehat{W}_{n-2}}{z}
\]
belongs to $\u(3)^{\C}_{n-2}$.  And since $\pi_0^\perp \in \u(3)^{\C}_2$,
$\widehat{W}_n$ is in $\gC_n$.

We now prove that (\ref{eq:champ_formel}) is also true along $\partial/\partial
\bar{z}$.  We follow here the same kind of arguments as in \cite{BFPP} slightly
simplified\footnote{essentially the simplifications occur because the
semi-simple term $B$ of \cite{BFPP} is here $ia\pi_0^\perp $ which is constant
for $d$, so we do not need to introduce a flat connection.}.  We want to show
that
\[
\frac{\partial Y_\lambda}{\partial \bar{z}} +
\left[A_\lambda\left(\frac{\partial }{\partial \bar{z}}\right), Y_\lambda\right]
= 0
\]
and for that purpose we rather consider the conjugate of the left hand side
$\zeta_\lambda = (1+W_\lambda) (\frac{\partial Y_\lambda}{\partial \bar{z}} +
[A_\lambda\left(\frac{\partial }{\partial \bar{z}}\right),
Y_\lambda])(1+W_\lambda)^{-1}$.  We then prove two facts
\begin{itemize}
\item $\zeta_\lambda$ takes its values in $V^\perp$: this follows from the
identity
\[
\zeta_\lambda = \left[\lambda^{-2}ia\pi_0^\perp , \frac{\partial
W_\lambda}{\partial \bar{z}} (1+W_\lambda)^{-1} -
(1+W_\lambda)A_\lambda\left(\frac{\partial }{\partial \bar{z}}\right)
(1+W_\lambda)^{-1}\right].
\]
Note that since $\zeta_\lambda$ is twisted the fact that $\zeta_\lambda\in
V^\perp$ implies also that $\zeta_\lambda$ is an odd function of $\lambda$
and so that
\begin{equation}\label{3.1.oddzeta}
\zeta_\lambda = \sum_{k=0}^\infty \widehat{\zeta}_{2k-1}\lambda^{2k-1}.
\end{equation}

\item the relation
\begin{equation}\label{3.1.dzeta}
{\partial \zeta_\lambda\over \partial z} = [\varphi_\lambda,\zeta_\lambda].
\end{equation}
Indeed $d+\ad A_\lambda$ is a flat connection and in particular
$\frac{\partial }{\partial z} +\ad A_\lambda(\frac{\partial }{\partial z})$
commutes with $\frac{\partial }{\partial \bar{z}} +A_\lambda(\frac{\partial
}{\partial \bar{z}})$.  Hence
\[
\left( \frac{\partial }{\partial z}+\ad A_\lambda\left(\frac{\partial
}{\partial z}\right) \right) \left( \frac{\partial }{\partial \bar{z}} +\ad
A_\lambda\left(\frac{\partial }{\partial \bar{z}}\right) \right) Y_\lambda
= 0
\]
i.e.
\[
\left( \frac{\partial }{\partial z}+\ad A_\lambda\left(\frac{\partial
}{\partial z}\right) \right)
\left((1+W_\lambda)^{-1}\zeta_\lambda(1+W_\lambda)\right) = 0
\]
Thus (\ref{3.1.dzeta}) follows from a computation which uses
$\varphi_\lambda = {\partial W_\lambda\over\partial z}(1+W_\lambda)^{-1} -
(1+W_\lambda) A_\lambda\left(\frac{\partial }{\partial
z}\right)(1+W_\lambda)^{-1}$.

\end{itemize}
Now assume by contradiction that $\zeta_\lambda\neq 0$: in view of
(\ref{3.1.oddzeta}) there exists an integer $k\in \mathbb{N}$ such that
$\widehat{\zeta}_{2k-1}\neq 0$ and $\widehat{\zeta}_{2k-3} = 0$.  By
substituting the Fourier decompositions in (\ref{3.1.dzeta}) and observing that
the Fourier series expansion of $\varphi_\lambda$ starts by $\lambda^{-2}
ia\pi_0^\perp$, we deduce that $0 = \partial \widehat{\zeta}_{2k-3}/\partial z =
[ia\pi_0^\perp , \widehat{\zeta}_{2k-1}]$; but $\ad \pi_0^\perp $ is invertible
on $V^\perp$ and hence $\widehat{\zeta}_{2k-1} = 0$.  So we get a contradiction.

\subsection{Polynomial Killing fields}
We now deduce the existence of a non-trivial polynomial Killing field.

A first easy consequence of the results of the previous section is that, for all
$n\in \N$ and for all polynomial of the form $P(\lambda)=a_n\lambda ^{-4n} +
a_{n-1}\lambda ^{-4(n-1)} + \cdots + a_0$, where $a_0,a_1,\cdots ,a_n\in \C$ and
$a_n\neq 0$, then $Z_\lambda:=P(\lambda)Y_\lambda$ is again formal Killing
field.  Moreover it is quasi-adapted (modulo the multiplicative factor
$a_n\lambda^{-4n}$), i.e.\,the lower degree terms are $a_n\lambda^{-4n}\left(ia
\lambda^{-2} \pi_0^\perp + \lambda^{-1} X + O(\lambda ^0)\right)$).  Let us
consider
\[
Z_\leq := \sum_{k=-2-4n}^0\widehat{Z}_k\lambda^k,\quad \hbox{and}\quad Z_>:=
\sum_{k=1}^\infty \widehat{Z}_k\lambda^k,
\]
so that $Z_\lambda= Z_\leq + Z_>$.  We study
\begin{equation}\label{dZ-}
R_\lambda:= dZ_\leq +[A_\lambda,Z_\leq].
\end{equation}
We first remark that $R_\lambda$ is necessarily of the form
$R_\lambda=\sum_{k=-4-4n}^2\widehat{R}_k\lambda^k$.  But because of $dZ_\lambda
+ [A_\lambda,Z_\lambda]= 0$, we also have
\begin{equation}\label{dZ+}
R_\lambda = -dZ_> - [A_\lambda,Z_>],
\end{equation}
which implies $R_\lambda=\sum_{k=-1}^\infty \widehat{R}_k\lambda^k$.  Hence
finally
\[
R_\lambda = \lambda^{-1}\widehat{R}_{-1} + \widehat{R}_{0} +
\lambda^{1}\widehat{R}_{1} + \lambda^{2}\widehat{R}_{2}.
\]
Each term $\widehat{R}_{k}$ can be evaluated through two different ways: by
using (\ref{dZ-}) or (\ref{dZ+}).  From (\ref{dZ-}) we obtain
\begin{equation}\label{dZ-z}
\left\{ \begin{array}{l} \widehat{R}_{-1}(\partial _z) = \partial _z
\widehat{Z}_{-1} + [A_{-1}(\partial _z), \widehat{Z}_0] + [A_{0}(\partial _z),
\widehat{Z}_{-1}] \\
\widehat{R}_{0}(\partial _z) = \partial _z \widehat{Z}_{0} + [A_{0}(\partial
_z), \widehat{Z}_0] \\
\widehat{R}_{1}(\partial _z) = 0\\
\widehat{R}_{2}(\partial _z) = 0 \end{array}\right.
\end{equation}
and
\begin{equation}\label{dZ-zbar}
\left\{ \begin{array}{l} \widehat{R}_{-1}(\partial _{\bar{z}}) = \partial
_{\bar{z}} \widehat{Z}_{-1} + [A_{0}(\partial_{\bar{z}}), \widehat{Z}_{-1}] +
[A_{1}(\partial_{\bar{z}}), \widehat{Z}_{-2}] + [A_{2}''(\partial _{\bar{z}}),
\widehat{Z}_{-3}] \\
\widehat{R}_{0}(\partial _{\bar{z}}) = \partial_{\bar{z}} \widehat{Z}_{0} +
[A_{0}(\partial _{\bar{z}}), \widehat{Z}_{0}] + [A_{1}(\partial _{\bar{z}}),
\widehat{Z}_{-1}] + [A_{2}''(\partial _{\bar{z}}), \widehat{Z}_{-2}] \\
\widehat{R}_{1}(\partial _{\bar{z}}) = [A_{1}(\partial _{\bar{z}}),
\widehat{Z}_{0}] + [A_{2}''(\partial _{\bar{z}}), \widehat{Z}_{-1}] \\
\widehat{R}_{2}(\partial _{\bar{z}}) = [A_{2}''(\partial _{\bar{z}}),
\widehat{Z}_{0}].  \end{array}\right.
\end{equation}
>From (\ref{dZ+}) we get
\begin{equation}\label{dZ+z}
\left\{ \begin{array}{l} \widehat{R}_{-1}(\partial _z) = - [A_{2}'(\partial _z),
\widehat{Z}_{1}] \\
\widehat{R}_{0}(\partial _z) = - [A_{2}'(\partial _z), \widehat{Z}_{2}] -
[A_{-1}(\partial _z), \widehat{Z}_{1}] \\
\widehat{R}_{1}(\partial _z) = - \partial _z\widehat{Z}_{1} - [A_{2}'(\partial
_z), \widehat{Z}_{3}] - [A_{-1}(\partial _z), \widehat{Z}_{2}] - [A_{0}(\partial
_z), \widehat{Z}_{1}] \\
\widehat{R}_{2}(\partial _z) = - \partial _z\widehat{Z}_{2} - [A_{2}'(\partial
_z), \widehat{Z}_{4}] - [A_{-1}(\partial _z), \widehat{Z}_{3}] - [A_{0}(\partial
_z), \widehat{Z}_{2}] \end{array}\right.
\end{equation}
and
\begin{equation}\label{dZ+zbar}
\left\{ \begin{array}{l} \widehat{R}_{-1}(\partial _{\bar{z}}) = 0 \\
\widehat{R}_{0}(\partial _{\bar{z}}) = 0 \\
\widehat{R}_{1}(\partial _{\bar{z}}) = - \partial _{\bar{z}} \widehat{Z}_{1} -
[A_{0}(\partial _{\bar{z}}), \widehat{Z}_{1}] \\
\widehat{R}_{2}(\partial _{\bar{z}}) = - \partial _{\bar{z}} \widehat{Z}_{2} -
[A_{0}(\partial _{\bar{z}}), \widehat{Z}_{2}] - [A_{1}(\partial _{\bar{z}}),
\widehat{Z}_{1}] .  \end{array}\right.
\end{equation}
Thus in order to obtain an expression of $R_\lambda$ which does depend only on
$\widehat{Z}_{-1}$ and $\widehat{Z}_{0}$, we exploit (\ref{dZ-z}) and the two
last equations in (\ref{dZ-zbar}).  But instead of using the two first equations
of (\ref{dZ-zbar}) we take the two first ones of (\ref{dZ+zbar}).  This gives us
\begin{equation}\label{Rz}
R_\lambda(\partial _z) = \lambda^{-1}\left( \partial _z\widehat{Z}_{-1} +
[A_{-1}(\partial _z),\widehat{Z}_0] + [A_{0}(\partial
_z),\widehat{Z}_{-1}]\right) + \left( \partial _z\widehat{Z}_{0} +
[A_{0}(\partial _z),\widehat{Z}_0] \right),
\end{equation}
\begin{equation}\label{Rzbar}
R_\lambda(\partial _{\bar{z}}) = \lambda \left( [A_{1}(\partial
_{\bar{z}}),\widehat{Z}_0] + [A_{2}''(\partial
_{\bar{z}}),\widehat{Z}_{-1}]\right) + \lambda^2[A_{2}''(\partial
_{\bar{z}}),\widehat{Z}_{0}].
\end{equation}
These relations will imply that $\widehat{Z}_{-1}$ and $\widehat{Z}_{0}$ satisfy
a second order elliptic equation.  In order to prove that we need to establish
another relation between $R_\lambda(\partial _z)$ and $R_\lambda(\partial
_{\bar{z}})$.  For that purpose recall that $dA_\lambda + A_\lambda\wedge
A_\lambda = 0$, which means that the connection $d+\ad A_\lambda$ has a
vanishing curvature.  In particular
\[
0= \left( d+\ad A_\lambda\right) \circ \left( d+\ad A_\lambda\right)Z_\leq =
dR_\lambda + [A_\lambda \wedge R_\lambda].
\]
This implies
\begin{equation}\label{fermeture}
{\partial R_\lambda(\partial _z)\over \partial \bar{z}} - {\partial
R_\lambda(\partial _{\bar{z}}) \over \partial z} = [A_\lambda(\partial
_z),R_\lambda(\partial _{\bar{z}})] - [A_\lambda(\partial
_{\bar{z}}),R_\lambda(\partial _z)].
\end{equation}
A substitution of (\ref{Rz}) and (\ref{Rzbar}) in (\ref{fermeture}) gives a
system of linear elliptic equations on $\widehat{Z}_{-1}$ and $\widehat{Z}_{0}$.
Since the space of solutions to this system which are periodic is finite
dimensional, it turns out that $\widehat{Z}_{-1}$ and $\widehat{Z}_{0}$ belong
to a finite dimensional vector space.  Hence relations (\ref{Rz}) and
(\ref{Rzbar}) force $R_\lambda(\partial _z)$ and $R_\lambda(\partial
_{\bar{z}})$ to stay in a finite dimensional vector space.

We can conclude: let us consider
\[
{\cal R}:= \{ R_\lambda / R_\lambda(\partial _z), R_\lambda(\partial
_{\bar{z}})\hbox{ are given by (\ref{Rz}) and (\ref{Rzbar}) and satisfy
(\ref{fermeture})}\}.
\]
It is a complex finite dimensional vector space.  Let us also denote by ${\cal
P}_n:=\{P(\lambda)=a_n\lambda ^{-4n} + a_{n-1}\lambda ^{-4(n-1)} + \cdots + a_0/
(a_0,\cdots , a_n)\in \C ^{n+1}\}$ and ${\cal P}_\infty :=\cup_{n\in \N}{\cal
P}_n$.

The linear map ${\cal P}_\infty \ni P(\lambda)\longmapsto dZ_\leq
+[A_\lambda,Z_\leq]$, where $Z_\leq = \left( P(\lambda)Y_\lambda\right)_\leq$
takes values in ${\cal R}$ and so has a finite rank, say $n$.  Then since
$\hbox{dim}_\mathbb{C} {\cal P}_n=n+1$, the map ${\cal P}_n \ni
P(\lambda)\longmapsto dZ_\leq +[A_\lambda,Z_\leq]$ has a non trivial kernel: let
$P(\lambda)=\sum_{k=0}^na_k\lambda^{-4k}$ be a non trivial polynomial in this
kernel.  Let $4p$ be the degree of $P$ in $\lambda^{-1}$, i.e.\,such that
$P(\lambda)=\sum_{k=0}^pa_k\lambda^{-4k}$ and $a_p\neq 0$.  Without loss of
generality we can assume that $a_p = 1$.  Then $\xi_\lambda:=
\left(P(\lambda)Y_\lambda\right)_\leq -
\left(P(\lambda)Y_\lambda\right)_\leq^\dagger$ is a solution of
(\ref{2.1.laxfn}) and (\ref{2.1.laxqfn2}).

\section{Homogeneous tori in $\mathbb{C} P^2$}

We describe here the simplest examples of Hamiltonian stationary Lagrangian tori
in $\mathbb{C}P^2$: the {{\em homogeneous}} Hamiltonian stationary Lagrangian
tori, i.e. immersions $u$ of $S^1 \times S^1$ into $\mathbb{C} P^2$ such that $u
( x + t, y ) = e^{t A} u ( x, y )$ and $u ( x, y + t ) = e^{t B} u ( x, y )$ for
some skew-Hermitian matrices $A$ and $B$.  Notice that $A$ and $B$ are only
defined up to addition with a multiple of $i \mathit{Id}$.  The simplest example
is the Clifford torus, namely the image by the Hopf map $\pi$ of the product
torus $\{ z = ( z^1, z^2, z^3 ) ; | z^1 | = | z^2 | = | z^3 | = 1 / \sqrt{3}
\}$.  This torus is minimal.  The main result states that all homogeneous
Hamiltonian stationary Lagrangian tori are similar to the Clifford torus.

\begin{theorem} \label{homogeneous}
Any homogeneous Hamiltonian stationary Lagrangian torus in $\mathbb{C} P^2$ is
the image by the Hopf map of some Cartesian product $r_1 S^1 \times r_2 S^1
\times r_3 S^1 = \{ z = ( z^1, z^2, z^3 ) ; | z^1 | = r_1, | z^2 | = r_2, | z^3
| = r_3 \}$ where $r_1^2 + r_2^2 + r_3^2 = 1$, up to $U(3)$ congruence.
Moreover, the torus is special Lagrangian if and only if $r_1 = r_2 = r_3 =
\sqrt{3}$.
\end{theorem}
\begin{proof}
Let us first see why $\pi ( T )$ is a Hamiltonian stationary Lagrangian torus in
$\mathbb{C} P^2$, where $T = r_1 S^1 \times r_2 S^1 \times r_3 S^1$.  Indeed it
suffices to show that $\pi ( T )$ admits a Legendrian preimage.  Let \[ f ( x, y
) : = \left( r_1 e^{i (( 1 - r_1^2 ) x - r_2^2 y )}, r_2 e^{i ( - r_1^2 x + ( 1
- r_2^2 ) y )}, r_3 e^{i ( - r_1^2 x - r_2^2 y )} \right) .
\]
Then the orbit under the Hopf action of the image of $f$ is exactly the 3-torus
$T$ above and $\pi \circ f$ is doubly periodic with periods $(2\pi,0)$ and
$(0,2\pi)$.  Note that this immersion is not conformal but there exists an
orthonormal Hermitian moving frame $(e_1,e_2)$ such that ${\partial f\over
\partial x} = r_1\sqrt{1-r_1^2}\, e_1$ and ${\partial f\over \partial y} =
\frac{r_2}{\sqrt{1-r_1^2}} (r_3 e_2 - r_1r_2 e_1)$.  And it is easy to check
that $f$ is Legendrian (and flat).  Its Lagrangian angle function is
\[
\beta ( x, y ) = x ( 1 - 3 r_1^2 ) + y ( 1 - 3 r_2^2 ) + \pi
\]
and since the metric is flat, $\beta$ is clearly harmonic, and constant if and
only if $r_1 = r_2 = r_3 = 1 / \sqrt{3}$.  Notice that many of these tori do not
lift up to $S^5$ as Legendrian tori (they do not close up).  Indeed the Maslov
class is not always an integer: for the implicit homology basis, $t \longmapsto
( 2 \pi t, 0 )$ and $t \longmapsto ( 0, 2 \pi t )$, it is $( 1 - 3 r_1^2, 1 - 3
r_2^2 )$.  However, if all $r_i^2$ are rational, the torus in $\mathbb{C} P^2$
possesses a Legendrian toric multiple cover.

Suppose now that $u : S^1 \times S^1 \longrightarrow \mathbb{C} P^2$ is a
homogeneous Lagrangian immersion.  According to our definition $u$ has a lift
$\hat{u}$ such that $\pi ( \hat{u} ( x + t, y )) = \pi ( e^{t A} \hat{u} ( x, y
))$ and $\pi ( \hat{u} ( x, y + t )) = \pi ( e^{t B} \hat{u} ( x, y ))$.  In
particular $\pi ( e^{x A} e^{y B} p ) = \pi ( e^{y B} e^{x A} p )$, for any $p
\in S^5$ in the image.  However the image is never contained in a complex
subspace of $\mathbb{C}^3$, hence $[ A, B ] \in i\mathbb{R} \mathit{Id}$.  Since
$[ A, B ]$ is traceless, $A$ and $B$ commute.

The obvious (non Legendrian) lift in $S^5$ is $( x, y ) \longmapsto e^{x A} e^{y
B} p$ where now $p = ( p_1, p_2, p_3 )$ is a fixed point mapped by the Hopf map
$\pi$ to $u ( 0, 0 )$.  A Legendrian lift $\hat{u}$ takes the following form:
$\hat{u} ( x, y ) = e^{i \theta ( x, y )} e^{x A} e^{y B} p$ for some function
$\theta$.  The horizontality condition implies $\langle ( i\frac{\partial
\theta}{\partial x} \mathit{Id} + A ) p, p \rangle_{\mathbb{C}^3} = 0$ so that
$\frac{\partial \theta}{\partial x} = i \frac{\langle A p, p
\rangle_{\mathbb{C}^3}}{| p |^2}$ is a constant.  The same holds in the $y$
direction so we can define the lift $\hat{u} ( x, y ) = e^{x \hat{A} + y
\hat{B}} p$ where $\hat{A} = A + i\frac{\partial \theta}{\partial x}
\mathit{Id}$ and $\hat{B} = B + i\frac{\partial \theta}{\partial y} \mathit{Id}$
are two commuting skew-symmetric matrices.  (Notice that $\hat{u}$ is only
defined on the universal cover $\mathbb{R}^2$.)  The base point $p$ depends of
course on the choice of origin and is only defined up to multiplication by a
complex unit number.  Nevertheless it plays an important role.

Consider now the metric induced by $\hat{u}$.  Due to homogeneity, it is a
constant metric on the $( x, y )$-plane.  By doing a simple change in variables,
we may as well assume that the metric is the standard plane metric, in other
words the immersion is isometric.  (Of course that will change the matrices
$\hat{A}$ and $\hat{B}$, but since they are replaced by some real linear
combination of themselves, the properties mentioned above still hold.)
Henceforth we suppose that $\hat{u}$ is an isometric homogeneous Legendrian
immersion of the plane.

Up to a unitary rotation in $\mathbb{C}^3$ we may suppose that $\hat{A}$ is
diagonal, and write $\hat{A} = i \,\hbox{diag} ( a_1, a_2, a_3 )$ with real
coefficients $a_1, a_2, a_3$.  We will now consider three cases and show that
only case (i) is possible.
\begin{enumerate}
\item Suppose $B = i \hbox{diag} ( b_1, b_2, b_3 )$ is diagonal.  Then the
surface lies inside the three torus $T = | p_1 | S^1 \times | p_2 | S^1
\times | p_3 | S^1$.  Necessarily it lifts $\pi ( T )$.  Isometry will
constrain the coefficients to be as above.

\item One and only one of the off-diagonal coefficients of $B$ is non zero.
We can assume it is $b_{12}$ up to permutation of the coordinates.
Commutation of $\hat{A}$ and $\hat{B}$ forces $a_1 = a_2$, while $a_3 \neq
a_1$, otherwise we would get a contradiction: $A$ cannot be a multiple of
$i\mathit{Id}$.  Let us first look at equations involving $A$.  The
immersion being isometric in $S^5$, $| p | = | A p | = 1$
\[
1 = | p_1 |^2 + | p_2 |^2 + | p_3 |^2 = a_1^2 ( | p_1 |^2 + | p_2 |^2 ) +
a_3^2 | p_3 |^2 \hspace{0.75em}
\]
but it is also Legendrian, so \[ \omega ( A p, p ) = \langle i A p, p
\rangle_{\mathbb{C}^3} = a_1 ( | p_1 |^2 + | p_2 |^2 ) + a_3 | p_3 |^2 =
0 .
\]
Hence \[ | p_1 |^2 + | p_2 |^2 = \frac{a_3}{a_3 - a_1} \hspace{0.25em},
\hspace{0.75em} | p_3 |^2 = - \frac{a_1}{a_3 - a_1} \hbox{ and } a_1 a_3
= - 1
\]
excluding thus $a_1 = 0$, and finally \[ | p_1 |^2 + | p_2 |^2 =
\frac{1}{1 + a_1^2} \hspace{0.75em}, \hspace{0.75em} | p_3 |^2 =
\frac{a_1^2}{1 + a_1^2} \hspace{0.75em} .
\]
Take now into account the Legendrian constraints on $B$: \[ B =
\left(\begin{array}{ccc} i b_1 & b_{12} & 0\\
- \overline{b_{12}} & i b_2 & 0\\
0 & 0 & b_3 \end{array}\right)
\]
\[ 0 = \langle B p, p \rangle_{\mathbb{C}^3} = i \left( \sum_1^3 b_j |
p_j |^2 + 2 \mathrm{Im} ( b_{12} \overline{p_1} p_2 ) \right) \] \[ 0 =
\langle B p, A p \rangle_{\mathbb{C}^3} = \sum_1^3 a_j b_j | p_j |^2 + 2
a_1 \mathrm{Im} ( b_{12} \overline{p_1} p_2 ).
\]
Uniting both, we deduce \[ a_1 b_1 | p_1 |^2 + a_1 b_2 | p_2 |^2 + a_1 b_3
| p_3 |^2 = a_1 b_1 | p_1 |^2 + a_1 b_2 | p_2 |^2 + a_3 b_3 | p_3 |^2.
\]
Since $| p_3 | \neq 0$ and $a_1 \neq a_3$, $b_3$ vanishes.  The remaining
equations are: \[ | p_1 |^2 + | p_2 |^2 = \frac{1}{1 + a_1^2}
\]
\begin{equation}
\label{bmatrix1} b_1 | p_1 |^2 + b_2 | p_2 |^2 + 2 \mathrm{Im} ( b_{12}
\overline{p_1} p_2 ) = 0
\end{equation}
\begin{equation}
\label{bmatrix2} ( b_1^2 + | b_{12} |^2 ) | p_1 |^2 + ( b_2^2 + |
b_{12}|^2 ) | p_2 |^2 + 2 ( b_1 + b_2 ) \mathrm{Im} ( b_{12}
\overline{p_1} p_2 ) = 1.
\end{equation}
Finally we infer a contradiction: indeed equation (\ref{bmatrix1}) amounts to
the existence of an isotropic vector $( p_1, p_2 )$ for the skew-hermitian
matrix $\left( \begin{array}{cc} i b_1 & b_{12}\\
- \overline{b_{12}} & i b_2 \end{array} \right)$, and that requires its
determinant $|b_{12}|^2 - b_1 b_2$ to vanish.  Plugging this into
(\ref{bmatrix2}), we obtain \[ ( b_1 + b_2 ) \big( b_1 | p_1 |^2 + b_2 |
p_2 |^2 + 2 \mathrm{Im} ( b_{12} \overline{p_1} p_2 ) \big) = 1,
\]
obviously contradicting (\ref{bmatrix1}).

\item If at least two off-diagonal coefficients of $B$ are non-zero, then $A =
ia_1 \mathit{Id}$.  But that contradicts $( A p | p ) = 0$.  So that case is
also excluded.
\end{enumerate}
\medskip

Notice that in the language of integrable systems, homogeneous tori correspond
to vacuum solutions and are of finite type for $p=0$.
\end{proof}

\affiliationone{ Fr\'ed\'eric H\'elein\\
Universit\'e Denis Diderot (Paris 7)\\
Institut de Math\'ematiques de Jussieu -- UMR 7586, Case 7012\\
2, place Jussieu\\
75251 Paris Cedex 05\\
France \email{helein@math.jussieu.fr}}
\affiliationtwo{ Pascal Romon\\
Universit\'e de Marne-la-Vall\'ee\\
5, bd Descartes, Champs-sur-Marne\\
77454 Marne-la-Vall\'ee Cedex 2 \\
France \email{romon@univ-mlv.fr}}

\end{document}